\pdfoutput=1
\RequirePackage{ifpdf}
\ifpdf 
\documentclass[pdftex]{sigma}
\else
\documentclass{sigma}
\fi

\numberwithin{equation}{section}
\newtheorem{Theorem}{Theorem}[section]
\newtheorem{Corollary}[Theorem]{Corollary}
\newtheorem{Lemma}[Theorem]{Lemma}
\newtheorem{Proposition}[Theorem]{Proposition}
\newtheorem{Conjecture}[Theorem]{Conjecture}
\newtheorem{Question}[Theorem]{Question}
{\theoremstyle{definition}
\newtheorem{Definition}[Theorem]{Definition}

\newtheorem{Example}[Theorem]{Example}
\newtheorem{Remark}[Theorem]{Remark}
\newtheorem{Remarks}[Theorem]{Remarks}
}
\DeclareMathOperator{\Area}{Area}
\DeclareMathOperator{\Comp}{Comp}
\DeclareMathOperator{\closure}{closure}
\DeclareMathOperator{\Det}{Det}
\DeclareMathOperator{\Diff}{Dif\/f}
\DeclareMathOperator{\Homeo}{Homeo}
\DeclareMathOperator{\Imop}{Im}
\DeclareMathOperator{\In}{In}
\DeclareMathOperator{\image}{image}
\DeclareMathOperator{\Image}{Image}
\DeclareMathOperator{\kernel}{kernel}
\DeclareMathOperator{\Loop}{Loop}
\DeclareMathOperator{\Out}{Out}
\DeclareMathOperator{\Prob}{Prob}
\DeclareMathOperator{\proj}{proj}
\DeclareMathOperator{\PSU}{PSU}
\DeclareMathOperator{\Reop}{Re}
\DeclareMathOperator{\Res}{Res}
\DeclareMathOperator{\Rot}{Rot}
\DeclareMathOperator{\spanop}{span}
\DeclareMathOperator{\weight}{weight}

\begin{document}

\newcommand{\arXivNumber}{1401.2675}

\allowdisplaybreaks

\renewcommand{\PaperNumber}{081}

\FirstPageHeading

\ShortArticleName{Werner's Measure on Self-Avoiding Loops and Welding}

\ArticleName{Werner's Measure on Self-Avoiding Loops\\
and Welding}

\Author{Angel CHAVEZ and Doug PICKRELL}
 \AuthorNameForHeading{A.~Chavez and D.~Pickrell}
\Address{Mathematics Department, University of Arizona, Tucson, AZ~85721, USA}
\Email{\href{mailto:achavez@math.arizona.edu}{achavez@math.arizona.edu},
\href{pickrell@math.arizona.edu}{pickrell@math.arizona.edu}} \URLaddress{\url{http://math.arizona.edu/~achavez/},
\url{http://math.arizona.edu/~pickrell/}}

\ArticleDates{Received February 18, 2014, in f\/inal form July 31, 2014; Published online August 04, 2014}

\Abstract{Werner's conformally invariant family of measures on self-avoiding loops on Riemann surfaces is determined~by
a~single measure $\mu_0$ on self-avoiding loops in $\mathbb C \setminus\{0\}$ which surround~$0$.
Our f\/irst major objective is to show that the measure $\mu_0$ is inf\/initesimally invariant with respect to conformal
vector f\/ields (essentially the Virasoro algebra of conformal f\/ield theory).
This makes essential use of classical variational formulas of Duren and Schif\/fer, which we recast in representation
theoretic terms for ef\/f\/icient computation.
We secondly show how these formulas can be used to calculate (in principle, and sometimes explicitly) quantities (such
as moments for coef\/f\/icients of univalent functions) associated to the conformal welding for a~self-avoiding loop.
This gives an alternate proof of the uniqueness of Werner's measure.
We also attempt to use these variational formulas to derive a~dif\/ferential equation for the (Laplace transform of) the
``diagonal distribution'' for the conformal \mbox{welding} associated to a~loop; this generalizes in a~suggestive way to
a~deformation of Werner's measure conjectured to exist by Kontsevich and Suhov (a basic inspiration for this paper).}


\Classification{60D05; 60B15; 17B68; 30C99}

\section{Introduction}

Given a~topological space~$S$, let $\Comp(S)$ denote the set of all compact subsets of~$S$ with the Vietoris
topology, and let
\begin{gather*}
\Loop(S):=\big\{\gamma\in \Comp(S): \gamma~\text{is homeomorphic to}~S^1\big\}
\end{gather*}
with the induced topology (see Appendix~\ref{appendixA}).
Suppose that for each Riemann surface~$S$, $\mu_S$~is a~positive Borel measure on $\Loop(S)$.
Following Werner, this family of measures is said to satisfy conformal restriction if for each conformal embedding
$S_1\to S_2$, the restriction of $\mu_{S_2}$ to $\Loop(S_1)$ equals $\mu_{S_1}$; the family is nontrivial if the
measure of the set
\begin{gather}
\label{nontrivial}
\big\{\gamma\in\Loop(\{0<\vert z\vert<A\})\setminus \Loop(\{\vert z\vert<a\}):\gamma~\text{surrounds}~0\big\}
\end{gather}
is f\/inite and positive, for some $0<a<A$.
In~\cite{W} Werner proved the following remarkable result.

\begin{Theorem}
There exists a~nontrivial family of measures $\{\mu_S\}$ on self-avoiding loops on Riemann surfaces which satisfies
conformal restriction.
This family is unique up to multiplication by an overall positive constant.
\end{Theorem}

\begin{Remark}
In the case of the plane, $S=\mathbb C$, this measure is conjectured to be the continuum limit of a~(properly weighted)
random self-avoiding circular walk (see Section~7.1 of~\cite{W}).
Kontsevich and Suhov have conjectured that there is a~deformation of Werner's (family of) measures to a~(family of)
measures having values in a~determinant line bundle $\Det^c$, where $c\le1$ is central charge; this deformation
is presumably a~continuum limit for other statistical mechanical models, and the determinant twist (when $c\ne0$) is
essential in understanding how the theory naturally extends to all Riemann surfaces (see~\cite{KS} and~\cite{BD}).
We refer the reader to the introduction of~\cite{W} and Section~6 of~\cite{KS} for further background and motivation,
and~\cite{Bauer} for a~connection with Schramm--Loewner evolution.
\end{Remark}

Below we will introduce a~normalization which uniquely determines Werner's family of measures
(see~\eqref{decomposition}).
We will assume this is in force from now on.

Essentially because any self-avoiding loop on a~Riemann surface is contained in an embedded annulus, the family
$\{\mu_S\}$ is (in principle) uniquely determined~by
\begin{gather*}
\mu_0:=\mu\vert_{\Loop^1(\mathbb C \setminus\{0\})}
\end{gather*}
the restriction of~$\mu$ to loops in the plane which surround~$0$.
The measure $\mu_0$ is determined, up to a~constant, by the following formula of Werner (see Proposition~3 of~\cite{W}):

\begin{Theorem}
Suppose that $0\in U\subset V$, where~$U$ and~$V$ are bounded simply connected domains in~$\mathbb C$.
Then
\begin{gather*}
\mu_0\big(\Loop^1(V\setminus\{0\})\setminus \Loop(U)\big)=c_{\rm W} \log(\vert\phi'(0)\vert),
\end{gather*}
where $\phi:(U,0)\to(V,0)$ is a~conformal isomorphism.
\end{Theorem}

We will refer to $c_{\rm W}$ as Werner's constant, which depends on the normalization~\eqref{decomposition}.
At the present time we can only say that $c_{\rm W}\ge1$ (see Section~\ref{Wernerconstant}).

Our purpose is to explore other possible explicit formulas for $\mu_0$, especially in terms of welding.
To put this in perspective, it is convenient to slightly digress and recall the ``fundamental theorem of Welding'', and
some associated terminology (we recommend~\cite{Bishop} as a~basic reference).

\begin{Theorem}
\label{QC}
Suppose that~$\sigma$ is a~quasisymmetric homeomorphism of $S^1$.
Then
\begin{gather*}
\sigma =l\circ ma\circ u,
\end{gather*}
where
\begin{gather*}
u=z\bigg(1+\sum\limits_{n\ge1}u_nz^n\bigg)
\end{gather*}
is a~univalent holomorphic function in the open unit disk~$\Delta$, with quasiconformal extension to $\mathbb C \cup
\{\infty\}$, $m\in S^1$ is a~rotation, $0<a\le1$ is a~dilation, the mapping inverse to~$l$,
\begin{gather*}
L(z)=z\bigg(1+\sum\limits_{m\ge1}b_mz^{-m}\bigg)
\end{gather*}
is a~univalent holomorphic function on the open unit disk about infinity, $\Delta^*$, with quasiconformal extension to
$\mathbb C \cup \{\infty\}$, and the compatibility condition
\begin{gather*}
mau\big(S^1\big)=L\big(S^1\big)
\end{gather*}
holds.
This factorization is unique.
\end{Theorem}

\begin{Definition}
\label{weldingdefn}
A~homeomorphism~$\sigma$ of $S^1$ has a~triangular factorization (or welding) if $\sigma =l\circ ma\circ u$, where
\begin{gather*}
u(z)=z\bigg(1+\sum\limits_{n\ge1}u_nz^n\bigg)
\end{gather*}
is a~holomorphic function in~$\Delta$ with a~continuous extension to a~homeomorphism on $D:=\closure(\Delta)$,
$m\in S^1$ is a~rotation, $0<a$ is a~dilation, the mapping inverse to~$l$,
\begin{gather*}
L(z)=z\bigg(1+\sum\limits_{m\ge1}b_mz^{-m}\bigg),
\end{gather*}
is a~holomorphic function on $\Delta^*$ with a~continuous extension to a~homeomorphism on
$D^*=\closure(\Delta^*)$, and the compatibility condition $mau\big(S^1\big)=L\big(S^1\big)$ holds.
\end{Definition}

Quasisymmetric homeomorphisms have unique triangular factorizations.
For less regular homeomorphisms, there are additional suf\/f\/icient conditions for the existence of weldings
(see~\cite{AJKS} and references, and~\cite{AMT2}), but there are many examples of homeomorphisms which do not admit
weldings, and weldings which are not unique (see~\cite{Bishop}).

Suppose that $\gamma\in\Loop^1(\mathbb C \setminus\{0\})$.
By the Jordan curve theorem the complement of~$\gamma$ in $\mathbb C \cup\{\infty\}$ has two connected components,
$U_{\pm}$, so that
\begin{gather*}
\mathbb C \cup \{\infty\}=U_+ \sqcup \gamma \sqcup U_-,
\end{gather*}
where $0\in U_+$ and $\infty \in U_-$.
There are based conformal isomorphisms
\begin{gather*}
\phi_+:\ (\Delta,0) \to (U_+,0),
\qquad
\phi_-: \ (\Delta^*,\infty) \to (U_-,\infty).
\end{gather*}
The map $\phi_-$ can be uniquely determined by normalizing the Laurent expansion in $\vert z\vert>1$ to be of the form
\begin{gather*}
\phi_-(z)=\rho_{\infty}(\gamma)L(z),
\qquad
L(z)=z\bigg(1+\sum\limits_{n\ge1}b_nz^{-n}\bigg),
\end{gather*}
where $\rho_{\infty}(\gamma)>0$ is the transf\/inite diameter (see Chapters~16 and~17 of~\cite{Hille} for numerous
formulas for $\rho_{\infty}$).
The map $\phi_+$ can be similarly uniquely determined by normalizing its Taylor expansion to be of the form
\begin{gather*}
\phi_+(z)=\rho_0(\gamma)u(z),
\qquad
u(z)=z\bigg(1+\sum\limits_{n\ge1}u_nz^{n}\bigg),
\end{gather*}
where $\rho_0(\gamma)>0$ is called the conformal radius with respect to~$0$.
By a~theorem of Carath\'{e}odory (see Theorem~17.5.3 of~\cite{Hille}), both $\phi_{\pm}$ extend uniquely to
homeomorphisms of the closures of their domain and target.
This implies that the restrictions $\phi_{\pm}:S^1 \to \gamma$ are topological isomorphisms.
Thus there is a~well-def\/ined welding map
\begin{gather}
\label{welding}
W:\ \Loop^1(\mathbb C \setminus\{0\}) \to \big\{\sigma\in \Homeo^+\big(S^1\big):\sigma=lau\big\} \times \mathbb R^+:\gamma
\mapsto (\sigma(\gamma),\rho_{\infty}(\gamma)),
\end{gather}
where
\begin{gather*}
\sigma(\gamma,z):=\phi_-^{-1}(\phi_+(z))=lau,
\qquad
a(\gamma)=\frac{\rho_0(\gamma)}{\rho_{\infty}(\gamma)}
\end{gather*}
and~$l$ is the inverse mapping for~$L$.

\begin{Remarks}\quad
\begin{enumerate}\itemsep=0pt
\item[(a)] To clarify~\eqref{welding}, the~$\sigma$ image of~$W$ is by def\/inition the set of homeomorphisms which admit
a~triangular factorization with rotation $m=1$.
\item[(b)] The map~$W$ is not $1-1$ because triangular factorization fails (in a~dramatic way) to be unique (the source of
nonuniqueness: there exist homeomorphisms of the $2$-sphere which are conformal of\/f of a~Jordan curve, and which are not
linear fractional transformations; see~\cite{Bishop}).
\end{enumerate}
\end{Remarks}

A lofty goal (not in sight) is to calculate, in some explicit way, the image measure $W_*\mu_0$, and to show that
$\mu_0$ can be recovered from this image.
As we will see in Section~\ref{Weldingmap}, conformal invariance implies that
\begin{gather}
\label{decomposition}
d(W_*\mu_0)(\sigma,\rho_{\infty})=d\nu_0(\sigma) \times \frac{d\rho_{\infty}}{\rho_{\infty}},
\end{gather}
where $\nu_0$ is an inversion invariant f\/inite measure, which we normalize to have unit mass.
This reduces the task of computing~$W_*\mu_0$ to computing the inversion invariant probability measure~$\nu_0$.

In this paper our f\/irst major objective is to show that the measure $\mu_0$ is inf\/initesimally invariant with respect to
conformal vector f\/ields, essentially the Virasoro algebra of conformal f\/ield theory.
This makes essential use of classical variational formulas of Duren and Schif\/fer~\cite{DS}, which we reformulate in
representation theoretic terms for ef\/f\/icient computation.
We secondly show how conformal invariance can be used to calculate integrals with respect to the measure~$\nu_0$
($\nu_0$~is not itself conformally invariant, so this is a~nontrivial step).
We thirdly show how these formulas can be used to calculate the joint moments for the coef\/f\/icients of~$u$.
Since these coef\/f\/icients are bounded, these moments (in principle) determine the joint distributions for the
coef\/f\/icients.
This yields an alternate proof of the uniqueness of Werner's measure.
This is also potentially interesting because a~suf\/f\/iciently explicit calculation of the individual moments for~$u_n$
could yield a~probabilistic proof of the Bieberbach conjecture/de Branges theorem (as pointed out by a~referee, one must
also show the measure $\nu_0$ has dense support in a~suitable sense).
Our current procedure (which we have implemented numerically) has the virtue that it in principle systematically
calculates all joint moments; it has the drawback that to obtain a~general moment for~$u_N$, it has to calculate on the
order of~$p(N)$ joint moments for all~$u_n$ with $n<N$, where $p(N)$ is the partition function (which grows very
rapidly).
In any event a~certain fraction of the moments turn out to have remarkably simple expressions; for example:

\begin{Theorem}
\label{2moment}
\[
\int \vert u_n\vert^2 d\nu_0=\frac1{n+1}.
\]
\end{Theorem}

The coef\/f\/icients of~$u$ are well-known to be functionally dependent in a~very complicated way (see Chapter~11
of~\cite{Duren}).
For this reason it seems unlikely that one could calculate~$\nu_0$ in an explicit way in terms of these coordinates.
For this reason it is important to consider other quantities (and coordinates) associated with the welding
homeomorphism~$\sigma$.
For various reasons (see Remarks~\ref{remarkbelow} below), it is of special interest to calculate the ``diagonal
distribution'', i.e.~the distribution for~$a$ in the triangular factorization $\sigma=lau$.

\begin{Conjecture}
If $\nu_0$ is normalized to be a~probability measure, then
\begin{gather*}
\nu_0(\{\sigma: \exp(-x)\le a(\sigma)\le1\})= \exp\left(-\frac{\beta_0}{x}\right),
\qquad
x>0
\end{gather*}
for some constant $\beta_0<\frac{5\pi^2}{4}$.
\end{Conjecture}

\begin{Remarks}\label{remarkbelow}\quad
\begin{enumerate}\itemsep=0pt
\item[(a)] The motivating idea is to show that the Laplace transform of the diagonal distribution for~$\nu_0$ satisf\/ies
a~dif\/ferential equation, using the inf\/initesimal conformal invariance of $\mu_0$.

\item[(b)] This conjecture is closely related to Proposition~18 in~\cite{W}, concerning the measure of the set of nontrivial
loops in a~f\/inite type annulus, for which there is an explicit conjecture due to Cardy (see
Section~\ref{diagonaldistribution}).

\item[(c)] There is a~natural generalization of this conjecture to the deformation of Werner's measure which is conjectured to
exist in~\cite{KS} (see Section~\ref{KSconjecture} and see~\cite{BD} for recent progress on this conjecture).
Our hope is that this extended conjecture might be useful in proving existence of this deformation.
\end{enumerate}
\end{Remarks}

To close this introduction, we mention one obvious coordinate which should be investigated.
For a~homeomorphism~$\sigma$ of $S^1$, write
\begin{gather*}
\sigma\big(e^{i\theta}\big)=e^{i\Sigma(\theta)},
\end{gather*}
where the lift~$\Sigma$ is a~homeomorphism of $\mathbb R$ satisfying $\Sigma(\theta+2\pi)=\Sigma(\theta)+2\pi$;~$\Sigma$
is determined modulo $2\pi \mathbb Z$.
The $\nu_0$ distribution for~$\sigma$ is completely determined by the distributional derivative,
\begin{gather*}
\frac1{2\pi} d\Sigma,
\end{gather*}
which we view as a~probability measure on $S^1$.

Verblunsky discovered a~remarkable parameterization of probability measures on $S^1$.
To state the gist of the result simply (following~\cite{simon}), let $\Prob'(S^1)$ denote the set of probability
measures which are nontrivial, in the sense that their support is not a~f\/inite set.

\begin{Theorem}
The following map induces a~bijective correspondence:
\begin{gather*}
\Prob'\big(S^1\big)\to \prod\limits_{n=0}^{\infty}\Delta:\omega \to (\alpha_n),
\end{gather*}
where $p_0=1$, $p_1=z-\overline{\alpha_1}$, \dots\
are the monic orthogonal polynomials with respect to~$\omega$, and $\alpha_n=-\overline{p_{n+1}(0)}$.
\end{Theorem}

It is very striking that the image of this correspondence is a~product space, i.e.~the $\alpha_n$ are functionally
independent, in sharp contrast to the coef\/f\/icients $u_n$.
This suggests the following naive

\begin{Question}
Are the Verblunsky coefficients
\begin{gather*}
(\alpha_n)\in \prod\limits_{n=0}^{\infty}\Delta
\end{gather*}
independent random variables with respect to $\nu_0$?
\end{Question}

We have basically failed in trying to investigate this question numerically.

\subsection{Outline of the paper}

In Section~\ref{Weldingmap} we prove some basic facts about the welding map~$W$.
In Section~\ref{variationalformulas} we recall some classical variational formulas of Duren and Schif\/fer.
In Section~\ref{stresstensor} we discuss the inf\/initesimal action from a~representation theoretic point of view, and we
recast the Duren--Schif\/fer formulas in terms of generating functions, using a~stress-energy tensor formulation common in
conformal f\/ield theory.
In Section~\ref{infinitesimalinvariance} we establish the version of inf\/initesimal conformal invariance of $\mu_0$
needed for our purposes.
In Section~\ref{moments} we apply this to compute moments of the coef\/f\/icients of~$u$, and to give an alternate proof of
the uniqueness of Werner's family of measures.
In Section~\ref{diagonaldistribution} we discuss the relation between the diagonal distribution conjecture and
Proposition~18 of~\cite{W}, and outline a~strategy for a~proof; we also brief\/ly indicate how the conjecture generalizes
to the deformation which is conjectured by Kontsevich and Suhov to exist in~\cite{KS}.

\subsection{Notations and conventions}

Given a~complex number~$z$, we often write $z^*$ for the complex conjugate, especially when~$z$ is represented~by
a~complicated expression.

Given a~Laurent expansion $f(z)=\sum f_nz^n$, we write $f^*(z)=\sum (f_{n})^* z^{-n}$ (for $z\in S^1$, $f^*(z)=f(z)^*$).
We also write $f_-(z)=\sum\limits_{n<0}f_nz^n$, $f_+(z)=\sum\limits_{n\ge0}f_nz^n$,
$f_{++}(z)=\sum\limits_{n>0}f_nz^n$, and $f_{-1}=\Res(f(z),z=0)$.

Throughout this paper, we view vector f\/ields on a~manifold as the Lie algebra of dif\/feomor\-phisms of the manifold; the
induced bracket is the negative of the usual bracket obtained by viewing vector f\/ields as derivations of functions on
the manifold.

\section{The welding map}
\label{Weldingmap}

In this section we consider the welding map~\eqref{welding}.

\begin{Proposition}\label{introlemma}\quad
\begin{enumerate}\itemsep=0pt
\item[$(a)$] The distributions for $\rho_0$ and $\rho_{\infty}$ are invariant with respect to dilation, i.e.~equivalent to Haar
measure for $\mathbb R^+$.

\item[$(b)$] $
d(W_*\mu_0)(\sigma,\rho_{\infty})=d\nu_0(\sigma) \times \frac{d\rho_{\infty}}{\rho_{\infty}},
$
where $\nu_0$ is a~finite measure $($which we will normalize to have unit mass$)$.

\item[$(c)$] The measure $d\nu_0(\sigma)$ is inversion invariant and invariant with respect to conjugation by $C:z\mapsto z^*$.

\item[$(d)$] The measure $d\nu_0(\sigma)$ is supported on~$\sigma$ having triangular factorization $\sigma=lau$, i.e.~$m=1$.

\item[$(e)$] For any $\gamma\in\Loop^1(\mathbb C \setminus\{0\})$,
\begin{gather*}
a(\sigma(\gamma))=\left(\frac{1-\sum\limits_{m=1}^{\infty}(m-1)\vert b_m\vert^2}{1+\sum\limits_{n=1}^{\infty}(n+1)\vert
u_n\vert^2}\right)^{1/2}\le1,
\end{gather*}
where~$u$ and~$L$ are written as in Definition~{\rm \ref{weldingdefn}}.

\item[$(f)$] The welding map is equivariant with respect to rotations in the sense that
\begin{gather*}
\sigma(\Rot(\theta)(\gamma))=\Rot(\theta)\circ \sigma(\gamma)\circ \Rot(\theta)^{-1}.
\end{gather*}
\end{enumerate}
\end{Proposition}

\begin{proof}
We f\/irst claim that
\begin{gather*}
\big\{\gamma\in \Loop^1(\mathbb C \setminus\{0\}): r<\rho_{\infty}(\gamma)<R\big\}\subset \Loop^1(\{\vert
z\vert<4R\})\setminus \Loop^1(\{r<\vert z\vert\}).
\end{gather*}
The inequality $r<\rho_{\infty}(\gamma)$ implies that~$\gamma$ cannot be contained in $\{\vert z\vert<r\}$.
In general $\rho_{\infty}(\gamma)=\rho_0(\frac1{\gamma})^{-1}$.
Thus if $\rho_{\infty}(\gamma)<R$, then $\frac1{R}<\rho_0(\frac1{\gamma})$.
The Koebe one-quarter theorem implies that $\frac1{\gamma}\subset \{\frac1{4R}<\vert z\vert\}$.
Thus~$\gamma$ is in the ball of radius~$4R$.
This proves the claim.

By conformal invariance and the nontriviality assumption of Werner, the set of loops (surrounding zero) with
$r<\rho_{\infty}<R$ has $\mu_0$ f\/inite measure, for any $r<R$.
This implies that there is a~essentially unique disintegration of $\mu_0$ of the form
\begin{gather*}
d\mu_0(\gamma)=\int_{\rho_{\infty}=0}^{\infty}d\mu_{\rho_{\infty}}(\sigma)d\omega(\rho_{\infty}),
\end{gather*}
where the f\/iber measures are probability measures.

The invariance of $\mu_0$ with respect to dilation, $\gamma \to \rho \gamma$, implies that the $\rho_{\infty}$
distribution~$\omega$ is also dilation invariant, i.e.~it is a~Haar measure for~$\mathbb R^+$.
The invariance of $\mu_0$ with respect to $z\to \frac1{z}$ implies that the same is true for $\rho_{0}$.
This proves (a).

Since $\mu_0$ is determined up to multiplication by a~constant, we can suppose that
\begin{gather*}
d\omega(\rho_{\infty})=\frac{d\rho_{\infty}}{\rho_{\infty}}.
\end{gather*}
The action by dilation transports one f\/iber to another.
Hence dilation invariance also implies that all the f\/iber measures are the same.
This implies that $W_*\mu_0$ is a~product measure, as claimed in part (b).

For part (c), we f\/irst use the invariance of $\mu_0$ with respect to $z\to \frac1{z^*}$, which maps~$\gamma$ to
$\frac1{\gamma^*}$:
\begin{gather*}
\phi_+\left(\frac1{\gamma^*}\right)(z)=\frac1{\phi_-(\gamma)\big(\frac1{z^*}\big)^*},
\qquad
\vert z\vert <1,
\\
\phi_-\left(\frac1{\gamma^*}\right)(z)=\frac1{\phi_+(\gamma)\big(\frac1{z^*}\big)^*},
\qquad
\vert z \vert>1,
\end{gather*}
and
\begin{gather*}
\phi_-\left(\frac1{\gamma^*}\right)^{-1}(w)=\frac1{\left(\phi_+(\gamma)^{-1}\big(\frac1{w^*}\big)\right)^*}.
\end{gather*}
Thus
\begin{gather*}
\sigma(\frac{1}{\gamma^*})(z)=\phi_-\left(\frac1{\gamma^*}\right)^{-1}\circ \phi_+\left(\frac1{\gamma^*}\right)(z)=
\phi_-\left(\frac1{\gamma^*}\right)^{-1}\left(\frac1{\phi_-(\gamma)\big(\frac1{z^*}\big)^*}\right)
\\
\hphantom{\sigma(\frac{1}{\gamma^*})(z)}{}
=\frac1{\left(\phi_+(\gamma)^{-1}(\phi_-(\gamma)\big(\frac1{z^*}\big)\right)^*}
=\frac1{\sigma^{-1}(\gamma)\big(\frac1{z^*}\big)^*}=\sigma(\gamma)^{-1}(z).
\end{gather*}
This implies the invariance of $\nu_0$ with respect to inversion.

The measure $\mu_0$ is also invariant with respect to $C:z\mapsto z^*$.
In this case
\begin{gather*}
\phi_{\pm}(\gamma^*)=C\circ \phi_{\pm}(\gamma)\circ C.
\end{gather*}
This implies that $\nu_0$ is invariant with respect to conjugation by~$C$.
This proves (c).

Part (d) is obvious.

For part (e) (essentially the well-known area theorem from the theory of univalent functions), the main point is that
\begin{gather*}
au(\Delta)=\mathbb C\setminus L(D^*).
\end{gather*}
For suf\/f\/iciently smooth~$\gamma$
\begin{gather*}
\Area(u(\Delta))=\frac{1}{2i}\int_{\gamma}d\bar{t}\wedge dt=\frac1{2i}\int_{S^1}\bar{u}du
=\pi\left(1+\sum\limits_{n=1}^{\infty}(n+1)\vert u_n\vert^2\right)
\end{gather*}
and
\begin{gather*}
\Area(\mathbb C\setminus L(\Delta^*))=\frac1{2i}\int_{\gamma}\bar{t}dt=\frac1{2i}\int_{S^1}\overline{L}dL
=\pi\left(1-\sum\limits_{m=1}^{\infty}(m-1)\vert b_m\vert^2\right).
\end{gather*}
By continuity of measure, these formulas hold for all~$\gamma$.
This implies part (e).

Part (f) follows from
\begin{gather*}
\phi_{\pm}(\Rot(\theta)(\gamma))=\Rot(\theta)\circ \phi_{\pm}\circ \Rot(\theta)^{-1}.\tag*{\qed}
\end{gather*}
\renewcommand{\qed}{}
\end{proof}

\begin{Remarks}\quad
\begin{enumerate}
\item[(a)] In connection with part (c), in general, if a~homeomorphism~$\sigma$ has a~triangular factorization $lmau$, then
$\sigma^{-1}$ has a~triangular factorization with
\begin{gather*}
u\big(\sigma^{-1}\big)(z)=\frac1{L\big(\frac1{z^*}\big)^*},
\qquad
l\big(\sigma^{-1}\big)(z)=\frac1{U\big(\frac1{z^*}\big)^*},
\\
m\big(\sigma^{-1}\big)=m(\sigma)^*,
\qquad
a\big(\sigma^{-1}\big)=a(\sigma).
\end{gather*}
In particular inversion stabilizes the set of~$\sigma$ having triangular factorization with $m=1$.

\item[(b)] In connection with part~(f), equivariance with respect to rotations, see Section~\ref{lackofequiv} below.
\end{enumerate}
\end{Remarks}

\subsection{Unresolved foundational issues}

Theorem~\ref{QC} implies that~$W$ induces a~bijection
\begin{gather*}
W: \ \text{QuasiCircles}^1(\mathbb C \setminus\{0\}) \leftrightarrow \big\{\sigma \in \text{QS}\big(S^1\big):\sigma=lau\big\}
\times\mathbb R^+,
\end{gather*}
where a~quasicircle is a~Jordan curve which admits a~parameterization by the restriction to $S^1$ of a~quasiconformal
homeomorphism of $\mathbb C\cup\{\infty\}$.

\begin{Conjecture}\quad
\begin{enumerate}\itemsep=0pt
\item[$(a)$] $\mu_0$ has measure zero on quasicircles.

\item[$(b)$] $W$ is $1-1$ on a~set of full $\mu_0$ measure.

\item[$(c)$] Almost surely with respect to $\nu_0$,~$\sigma$ has a~unique triangular factorization $($with $m=1)$.
\end{enumerate}
\end{Conjecture}

This kind of issue is addressed in~\cite{AMT2} and~\cite{AJKS}.

In this paper we will need to avoid these unresolved issues.
In particular, because we do not know that~$u$ is determined by~$\sigma$ (in an almost sure sense), in the remainder of
this paper, we will implicitly view $\nu_0$ as a~measure on $\{u\}$.
Thus in place of (b) of Proposition~\ref{introlemma}, we will use the following decomposition, which is proved in exactly the same way.

\begin{Proposition}\label{technical}
\begin{gather*}
d\mu_0(u,\rho_{\infty})=d\nu_0(u) \times \frac{d\rho_{\infty}}{\rho_{\infty}},
\end{gather*}
where $\nu_0$ is a~finite measure $($which we will normalize to have unit mass$)$.
\end{Proposition}

\section{Variational formulas}
\label{variationalformulas}

Since the measure $\mu_0$ has a~local form of conformal invariance, it is natural to suspect that there are senses in
which the measure is inf\/initesimally conformally invariant.
For this reason we need to consider how $\phi_{\pm}$ vary when the curve~$\gamma$ is varied by a~local deformation
$z\mapsto z+\epsilon v(z)$, where $v(z)$ is holomorphic in $\mathbb C\setminus \{0\}$.
This deformation corresponds to a~real vector f\/ield
\begin{gather*}
\overset{\rightarrow}{v}=v_1\frac{\partial}{\partial x}+v_2\frac{\partial}{\partial y},
\end{gather*}
where $v=v_1+iv_2$.
Let $\overset{\longrightarrow}{\mathcal W}$ denote the real Lie algebra of all such vector f\/ields, where $v(z)$ has
a~f\/inite Laurent expansion.

For technical reasons, we distinguish $\overset{\longrightarrow}{\mathcal W}$ from the Witt algebra $\mathcal W$, which
consists of holomorphic vector f\/ields $v(z)\frac{\partial}{\partial z}$, where again $v(z)$ has a~f\/inite Laurent
expansion.
The Witt algebra is a~complex Lie algebra.
It is spanned over $\mathbb C$ by the vector f\/ields
\begin{gather*}
L_n=-z^{n+1}\frac{\partial}{\partial z},
\qquad
n\in\mathbb Z
\end{gather*}
with bracket
\begin{gather*}
[L_n,L_m]=(m-n)L_{n+m}.
\end{gather*}
$W^*$ consists of antiholomorphic vector f\/ields.
It is spanned by $\{\overline{L}_n:n\in \mathbb Z\}$.

The precise relationship between $\overset{\longrightarrow}{\mathcal W}$ and $\mathcal W$ is that there is a~real
embedding
\begin{gather}
\label{realembedding}
\overset{\longrightarrow}{\mathcal W} \to \mathcal W \oplus \mathcal W^*: \
\overset{\rightarrow}{v}=v_1\frac{\partial}{\partial x}+v_2\frac{\partial}{\partial y} \to v(z)\frac{\partial}{\partial
z}+v^*(z)\frac{\partial}{\partial \overline{z}}.
\end{gather}
Loosely speaking, $\overset{\longrightarrow}{\mathcal W}$ is the Witt algebra considered as a~real Lie algebra (see~\cite[p.~115]{DMS}).
The reason for maintaining a~distinction is that the variational formulas below will naturally def\/ine a~real
representation of the real Lie algebra $\overset{\longrightarrow}{\mathcal W}$.
However it is convenient to express this representation in terms of an associated complex representation of~$\mathcal W$.
We will write the map~\eqref{realembedding} as
\begin{gather*}
\overset{\rightarrow}{L} \leftrightarrow (L,\overline{L}).
\end{gather*}
In particular
\begin{gather*}
\overset{\rightarrow}{L}_n \leftrightarrow (L_n,\overline{L}_n)
\qquad
\text{and}
\qquad
\overset{\longrightarrow}{iL}_n
\leftrightarrow(iL_n,\overline{iL}_n)=(iL_n,-i(\overline{L}_n)).
\end{gather*}

\subsection{Variational formulas, I}

$\overset{\longrightarrow}{\mathcal W}$ can be viewed as a~Lie algebra of vector f\/ields on $\Loop^1(\mathbb C\setminus\{0\})$,
where by def\/inition a~vector f\/ield on a~self-avoiding loop is simply a~$\mathbb R^2$-valued vector
f\/ield along the loop (the degree of smoothness of a~loop is not relevant here).
In particular
\begin{gather*}
\overset{\rightarrow}{L}_n\big\vert_{\gamma}=\frac{d}{dt}\exp \big(t\overset{\rightarrow}{L}_n\big)(\gamma)\vert_{t=0}
\end{gather*}
and similarly for $\overset{\longrightarrow}{iL}_n$.
The corresponding actions on a~function of~$\gamma$ are given~by{\samepage
\begin{gather*}
\overset{\rightarrow}{L}_n\cdot
F(\gamma)=\frac{d}{dt}\big\vert_{t=0}F\big(\exp \big({-}t\overset{\rightarrow}{L}_n\big)(\gamma)\big)=\frac{d}{dt}\big\vert_{t=0}F\big(\gamma+t\gamma^{n+1}\big),
\end{gather*}
where in the last line we have implicitly chosen a~parameterization for~$\gamma$, and similarly for
$\overset{\longrightarrow}{iL}_n$.}

When $n\ge -1$, $L_n=-z^{n+1}\frac{\partial}{\partial z}$ is regular at $z=0$.
In this case it is very easy to f\/ind the variations of $\phi_+$ with respect to $\overset{\rightarrow}{L}_n$ and
$\overset{\longrightarrow}{iL}_n$.

\begin{Proposition}\quad
\begin{enumerate}\itemsep=0pt
\item[$(a)$] For $n\ge0$, $\overset{\rightarrow}{L}_n\phi_+=\phi_+^{n+1}$.
In particular
\begin{gather*}
\overset{\rightarrow}{L}_0\rho_0=\rho_0,
\qquad
\text{and}
\qquad
\overset{\rightarrow}{L}_0 u_k=0,
\qquad
k\ge1,
\end{gather*}
and for $n>0$
\begin{gather*}
\overset{\rightarrow}{L}_n\rho_0=0,
\qquad
\overset{\rightarrow}{L}_n u_k=0
\qquad
k<n,
\qquad
\overset{\rightarrow}{L}_n u_n=\rho_0^n,
\end{gather*}
and for $n<k$
\begin{gather*}
\overset{\rightarrow}{L}_n u_{k}=\rho_0^n p_{k-n}^{(n+1)}(u_1,\dots,u_{k-n}),
\end{gather*}
where
\begin{gather*}
\big(1+u_1z+u_2z^2+\cdots\big)^{n+1}=\sum\limits_{l=0}^{\infty}p_l^{(n+1)}(u_1,\dots,u_l)z^l.
\end{gather*}

\item[$(b)$] $\overset{\rightarrow}{L}_{-1}\phi_+=1+u'(z)\big({-}1+(u_1-u_1^*)z+z^2\big)$.

\item[$(c)$] For $n>0$, $\big(\overset{\longrightarrow}{iL}_n\big)\phi_+=i\phi_+^{n+1}$, and
$\big(\overset{\longrightarrow}{iL}_0\big)\phi_+=i(\phi_+-z\phi_+')$.

\item[$(d)$] $(\overset{\longrightarrow}{iL}_{-1})\phi_+=i\big(1+u'(z)\big({-}1+(u_1+u_1^*)z-z^2\big)\big)$.
\end{enumerate}
\end{Proposition}

\begin{proof}
$\overset{\rightarrow}{L}_0$ is inf\/initesimal dilation.
In this case the formulas in part (a) are obvious, because $\sigma(\gamma)$ is unchanged when~$\gamma$ is dilated.

For $n\ge -1$ and~$\epsilon$ suf\/f\/iciently small, a~uniformization for the region inside $\gamma+\epsilon\gamma^{n+1}$ is
the composition
\begin{gather*}
\phi_++\epsilon (\phi_+)^{n+1}.
\end{gather*}
This uniformization has to be composed with a~linear fractional transformation to obtain the correct normalization.
Consequently
\begin{gather}
\label{normalmapping}
\phi_+\big(\gamma+\epsilon \gamma^{n+1}\big)(z)=\big(\phi_++\epsilon
(\phi_+)^{n+1}\big)\left(\frac{\lambda(\epsilon)z+\omega(\epsilon)}{1+\omega(\epsilon)^*\lambda(\epsilon)z}\right),
\end{gather}
where~$\lambda$ (having unit norm) and~$\omega$ are determined by the conditions that this uniformization vanishes at
$z=0$ and has positive derivative at $z=0$.

Suppose $n\ge0$.
In this case the linear fractional transformation is the identity for all~$\epsilon$.
This implies part~(a).

Part (b), when $n+1=0$, is slightly more involved.
In this case
\begin{gather}
\big(\overset{\rightarrow}{L}_{-1}\phi_+\big)(z)
=\frac{d}{d\epsilon}\big\vert_{\epsilon=0}
\left(\phi_+\left(\frac{\lambda(\epsilon)z+\omega(\epsilon)}{1+\omega(\epsilon)^*\lambda(\epsilon)z}\right)+\epsilon\right)
\nonumber
\\
\phantom{\big(\overset{\rightarrow}{L}_{-1}\phi_+\big)(z)}
\label{calculation1}
=\phi_+'(z)\big(\dot{\lambda}(0)z+\dot{\omega}(0)-z\dot{\omega}(0)^* z\big)+1.
\end{gather}
To calculate the derivatives at zero, we use the normalizations for the mapping~\eqref{normalmapping}.
Because~$0$ is mapped to zero,
\begin{gather*}
\epsilon+\phi_+(\omega(\epsilon))=0.
\end{gather*}
This implies
\begin{gather*}
\omega(\epsilon)=\phi_+^{-1}(-\epsilon)
\qquad
\text{and}
\qquad
\dot{\omega}(0)=-\rho_0^{-1}.
\end{gather*}
Secondly the derivative of the map~\eqref{normalmapping} at $z=0$ must be positive.
Thus
\begin{gather*}
\phi_+'(\omega(\epsilon))\lambda(\epsilon)\big(1-\vert \omega(\epsilon)\vert^2\big)>0
\end{gather*}
and (because~$\lambda$ has unit norm)
\begin{gather*}
\lambda(\epsilon)^{-1}=\exp (i\Imop(\log (\phi_+'(\omega(\epsilon))))).
\end{gather*}
This implies
\begin{gather*}
\dot{\lambda}(0)
=-i\Imop\left(\frac{\phi_+''(\omega(0))\dot{\omega}(0)}{\phi_+'(\omega(0))}\right)=-2i\Imop\left(\frac{u_1}{\rho_0}\right).
\end{gather*}
Plugging these derivatives into~\eqref{calculation1} yields
\begin{gather*}
\overset{\rightarrow}{L}_{-1}\phi_+=1+\rho_0\big(1+2u_1 z+3u_2z^2+\cdots\big)\left(\frac{u_1-u_1^*}{\rho_0}z-\frac1{\rho_0}+\frac1{\rho_0}z^2\right).
\end{gather*}
This implies part (b).

Parts (c) and (d) are similar.
For part (c), when $n=0$, note that
\begin{gather*}
\big(\exp \big(i\theta \overset{\rightarrow}{L}_0\big)\phi_+\big)(z)=e^{i\theta}\phi_+\big(e^{-i\theta}z\big),
\end{gather*}
so that
\begin{gather*}
\exp \big(i\theta \overset{\rightarrow}{L}_0\big)\rho_0=\rho_0  \qquad\text{and}\qquad \exp \big(i\theta
\overset{\rightarrow}{L}_0\big)u_k=e^{-ik\theta}u_k.
\end{gather*}

For part (d), when $n+1=0$,
\begin{gather}
\label{normalmapping2+}
\phi_+(\gamma+i\epsilon)(z)=
\phi_+\left(\frac{\lambda(\epsilon)z+\omega(\epsilon)}{1+\omega(\epsilon)^*\lambda(\epsilon)z}\right)+i\epsilon
\end{gather}
and
\begin{gather}
\big(\overset{\longrightarrow}{iL}_{-1}\phi_+\big)(z)
=\frac{d}{d\epsilon}\big\vert_{\epsilon=0}\left(\phi_+\left(\frac{\lambda(\epsilon)z+\omega(\epsilon)}{1+\omega(\epsilon)^*\lambda(\epsilon)z}
\right)+i\epsilon\right)
\nonumber
\\
\phantom{\big(\overset{\longrightarrow}{iL}_{-1}\phi_+\big)(z)}
\label{calculation2}
=\phi_+'(z)\big(\dot{\lambda}(0)z+\dot{\omega}(0)-z\dot{\omega}(0)^* z\big)+i.
\end{gather}
To calculate the derivatives at zero, we use the normalizations for the mapping~\eqref{normalmapping2+}.
Because~$0$ is mapped to zero,
\begin{gather*}
i\epsilon+\phi_+(\omega(\epsilon))=0.
\end{gather*}
This implies
\begin{gather*}
\omega(\epsilon)=\phi_+^{-1}(-i\epsilon)
\qquad
\text{and}
\qquad
\dot{\omega}(0)=-i\rho_0^{-1}.
\end{gather*}
Secondly the derivative of the map~\eqref{normalmapping2+} at $z=0$ must be positive.
Thus
\begin{gather*}
\phi_+'(\omega(\epsilon))\lambda(\epsilon)\big(1-\vert \omega(\epsilon)\vert^2\big)>0
\end{gather*}
and (because~$\lambda$ has unit norm)
\begin{gather*}
\lambda(\epsilon)^{-1}=\exp\big(i\Imop(\log(\phi_+'(\omega(\epsilon))))\big).
\end{gather*}
This implies
\begin{gather*}
\dot{\lambda}(0)
=-i\Imop\left(\frac{\phi_+''(\omega(0))\dot{\omega}(0)}{\phi_+'(\omega(0))}\right)=2i\Imop\left(\frac{iu_1}{\rho_0}\right)
=i\frac{u_1+u_1^*}{\rho_0}.
\end{gather*}
Plugging these derivatives into~\eqref{calculation2} yields
\begin{gather*}
\overset{\longrightarrow}{iL}_{-1}\phi_+=i+\rho_0\big(1+2u_1 z+3u_2
z^2+\cdots\big)\left(\frac{i(u_1+u_1^*)}{\rho_0}z-\frac{i}{\rho_0}-\frac{i}{\rho_0}z^2\right).
\end{gather*}
This implies part (d).
\end{proof}

\subsection[(Lack of) Equivariance for~$W$]{(Lack of) Equivariance for~$\boldsymbol{W}$}
\label{lackofequiv}

We have already observed that the welding map is equivariant with respect to the actions of rotation of loops and
conjugation of homeomorphisms; see (f) of Proposition~\ref{introlemma}.

Given $\vert w\vert<1$, def\/ine $\phi_1(w)\in \PSU(1,1)$ (viewed as the group of automorphisms of the Riemann sphere which stabilize the circle)~by
\begin{gather*}
\phi_1(w;z)=\frac{z+\bar{w}}{1+w z}.
\end{gather*}

\begin{Proposition}
Suppose that $\gamma\in \Loop^1(\mathbb C\setminus\{0\})$ such that $\phi_1(\epsilon,\gamma)\in\Loop^1(\mathbb C\setminus\{0\})$.
Then to first order in~$\epsilon$
\begin{enumerate}\itemsep=0pt
\item[$(a)$] $\phi_+(\phi_1(\epsilon,\gamma))=\phi_1(\epsilon)\circ \phi_+(\gamma)\circ \phi_1(-\epsilon/\rho_{0})\circ
\exp(2i\Imop(u_1\bar{\epsilon})/\rho_0)$,

\item[$(b)$]
$\phi_-(\phi_1(\epsilon,\gamma))=\phi_1(\epsilon)\circ \phi_-(\gamma)\circ \phi_1(-\rho_{\infty}\epsilon)\circ
\exp(-2i\rho_{\infty}\Imop(b_1\epsilon))$,

\item[$(c)$] $\sigma(\phi_1(\epsilon,\gamma))= \exp(-2i\rho_{\infty}\Imop(b_1\epsilon))\circ\phi_1(\rho_{\infty}\epsilon)\circ
\sigma(\gamma)\circ \phi_1(-\epsilon/\rho_{0})\circ \exp(2i\Imop(u_1\bar{\epsilon})/\rho_0)$.
\end{enumerate}
\end{Proposition}

\begin{Remark}
The formula in (c) illustrates how the welding map is trying (with limited success) to intertwine the action of
$\PSU(1,1)$ on loops with its action by conjugation on the welding homeomorphism.
\end{Remark}

\begin{proof}
A~uniformization for the region inside $\gamma(\epsilon)$ is the composition
\begin{gather*}
\phi_1(\epsilon,\phi_+(\gamma)(z)).
\end{gather*}
This uniformization has to be precomposed with a~linear fractional transformation to obtain the correct normalization.
Consequently
\begin{gather*}
\phi_+(\gamma(\epsilon))= \phi_1(\epsilon)\circ\phi_+(\gamma)\circ\Phi_1(\epsilon),
\end{gather*}
where
\begin{gather*}
\Phi_1(\epsilon,z)=\frac{\lambda(\epsilon)z+\overline{\omega}(\epsilon)} {1+\omega(\epsilon)\lambda(\epsilon) z}
\end{gather*}
and~$\lambda$ (having unit norm) and~$\omega$ are determined by the conditions that this uniformization vanishes at
$z=0$ and has positive derivative at $z=0$.

The f\/irst condition implies
\begin{gather*}
\Phi_1(\epsilon,0)=\overline{\omega}(\epsilon)=\phi_+^{-1}(-\epsilon)=-\frac{\epsilon}{\rho_0}+O\big(\epsilon^2\big),
\end{gather*}
in particular $\omega'(0)=-\rho_0^{-1}$.
Note that for $\Phi_1$ to exist, $-\epsilon$ must be in $U_+$.
The second condition
\begin{gather*}
\phi_1(\epsilon)'[\phi_+\circ\Phi_1(0)]\phi_+'[\Phi_1(0)]\Phi_1'(0)>0
\end{gather*}
is equivalent to
\begin{gather*}
\phi_1(\epsilon)'[\phi_+(\overline{\omega}(\epsilon)]\phi_+'[\overline{\omega}(\epsilon)]\lambda \big(1-\vert
\omega\vert^2\big)>0
\end{gather*}
or
\begin{gather*}
\overline{\lambda}(\epsilon)=\frac{\phi_1(\epsilon)'[\phi_+(\overline{\omega}(\epsilon)]\phi_+'[\overline{\omega}(\epsilon)]}
{\vert\phi_1(\epsilon)'[\phi_+(\overline{\omega}(\epsilon)]\phi_+'[\overline{\omega}(\epsilon)] \vert}.
\end{gather*}
Use
\begin{gather*}
\phi_1(\epsilon)'(z)=\frac{1-\vert \epsilon \vert^2}{(1+\epsilon z)^2},
\\
\phi_+(\overline{\omega}(\epsilon))=\rho_0\overline{\omega}(\epsilon)+O\big(\epsilon^2\big)=-\epsilon+O\big(\epsilon^2\big),
\\
\phi_1(\epsilon)'(\phi_+(\overline{\omega}(\epsilon))= \frac{1-\vert \epsilon \vert^2}{\big(1+\epsilon
\big({-}\frac{\epsilon}{\rho_0}+O\big(\epsilon^2\big)\big)\big)^2}=1+O\big(\epsilon^2\big),
\\
\phi_+'[\overline{\omega}(\epsilon)]=\rho_0+2\rho_0u_1 \overline{\omega}(\epsilon)= \rho_0-2u_1 \epsilon+O\big(\epsilon^2\big).
\end{gather*}

Putting everything together
\begin{gather*}
\overline{\lambda}(\epsilon)=
\frac{(1+\epsilon^2+\cdots )(\rho_0-2u_1\epsilon+\cdots)}{\vert(1+\epsilon^2+\cdots)(\rho_0-2u_1\epsilon+\cdots)\vert}
=1-2\frac{(u_1-\overline{u_1})\epsilon}{\rho_0}+O\big(\epsilon^2\big).
\end{gather*}
This implies the formula in (a).

In a~similar way
\begin{gather*}
\phi_-(\phi_1(\epsilon;\gamma))=\phi(\epsilon)\circ\phi_-(\gamma)\Psi_1
\end{gather*}
and one precedes as before.
This leads to (b) and (c).
\end{proof}

\subsection{Variational formulas, II}

It is far more dif\/f\/icult to calculate $\overset{\rightarrow}{L}_{-n}\phi_+$ for $n>1$.
In this case $z^{-n+1}$ is regular at $z=\infty$.
This is the situation considered in~\cite{DS}, with slight modif\/ications.
The following statement is essentially equation (17) in~\cite{DS}.

\begin{Proposition}
Suppose that $n\ge1$.
Then
\begin{enumerate}\itemsep=0pt
\item[$(a)$] $\overset{\rightarrow}{L}_{-n}\rho_0=\rho_0^{-n+1}\Reop(P_n(u_1,\dots,u_n))$, where
\begin{gather*}
P_n(u_1,\dots,u_n)=\Res\left(\left(\frac{U'(t)}{U(t)}\right)^2t^{-n+1},t=0\right)
\end{gather*}
$($as always,~$U$ is the mapping inverse to~$u)$.
If $\deg (u_j)=j$, then $P_n$ is a~homogeneous polynomial of degree~$n$.

\item[$(b)$] For $k\ge1$
\begin{gather*}
\overset{\rightarrow}{L}_{-n}u_k=\rho_0^{-n}\left(ku_k\Reop(B_0)+\sum\limits_{m=1}^{k}(k+1-m)u_{k-m}(B_{-m}+\overline{B_m})\right),
\end{gather*}
where
\begin{gather*}
B_m=\Res\left(\left(\frac{U'(t)}{U(t)}\right)^2 U(t)^m t^{-n+1},t=0\right)=\Res\left(\frac{u(z)^{-n+1}}{u'(z)}z^{m-2},z=0\right).
\end{gather*}

\item[$(c)$] $\overset{\longrightarrow}{iL}_{-n}\rho_0=-\rho_0^{-n+1}\Imop(P_n(u_1,\dots,u_n))$.

\item[$(d)$] For $k\ge1$
\begin{gather*}
\overset{\longrightarrow}{iL}_{-n}u_k=-\rho_0^{-n}\left(ku_k\Imop(B_0)+
\sum\limits_{m=1}^{k}(k+1-m)u_{k-m}i(-B_{-m}+\overline{B_m})\right).
\end{gather*}
\end{enumerate}
\end{Proposition}

\begin{Remarks}
(a) It is natural to restate the relationship between the $P_n$ and~$U$ in terms of quadratic dif\/ferentials
\begin{gather*}
(\partial \log(U(t)))^2=\sum\limits_{n=0}^{\infty} P_{n}(u)t^{n}\left(\frac{dt}{t}\right)^2,
\end{gather*}
where $t=u(z)$, $z=U(t)$. In Section~\ref{diagonaldistribution} it will be convenient to rewrite this as
\begin{gather*}
\big(\partial \log\big(\phi_+^{-1}(t)\big)\big)^2=\sum\limits_{n=0}^{\infty} \rho_0^{-n}P_{n}(u)t^{n}\left(\frac{dt}{t}\right)^2
\end{gather*}
and to set $P_n(\phi_+)=\rho_0^{-n}P_{n}(u)$.
Hopefully this will not cause any confusion.

(b) Similarly the residue formula for $B_m$ is naturally understood as the integral over~$\gamma$ of the natural pairing
of the holomorphic vector f\/ield $-v(t)\frac{d}{dt}$ and the holomorphic quadratic dif\/ferential $U(t)^m(\partial
\log(U(t)))^2$.
\end{Remarks}

For later reference we note some elementary properties of the polynomials $P_n$.

\begin{Proposition}\quad
\begin{enumerate}\itemsep=0pt
\item[$(a)$] $P_n(u)$ is a~homogeneous polynomial in $u_1,\dots,u_n$ of degree~$n$, where $\deg(u_j)=j$, with integer coefficients.

\item[$(b)$] $P_n(u)=-2nu_n+\text{terms involving}~u_1,\dots,u_{n-1}$.

\item[$(c)$] $u_n$ is a~homogeneous polynomial in $P_1,\dots,P_n$ of degree~$n$, where $\deg(P_j)=j$, with rational coefficients.

\item[$(d)$] $u_n=-\frac{1}{2n}P_n+ \text{terms involving}~P_1,\dots,P_{n-1}$.
\end{enumerate}
\end{Proposition}

Thus $\mathbb Z[P_1,\dots,P_n]\subset \mathbb Z[u_1,\dots,u_n]$ is a~proper inclusion, but over $\mathbb Q$ they are the
same.

At this point we have formulas for the action of the real Witt algebra on the coef\/f\/icients of~$\phi_+$.
If we write
\begin{gather*}
\frac1{\phi_-(\frac1{w})}=\frac1{\rho_{\infty}}w\bigg(1+\sum\limits_{n\ge1}l_n w^{n}\bigg),
\end{gather*}
where $w=\frac1{z}$ is the standard coordinate at inf\/inity, then we can also write down formulas for the action of the
real Witt algebra on the coef\/f\/icients of~$\phi_-$.
We will postpone this until the next section.

\section{Reformulation of the variational formulas}
\label{stresstensor}

\subsection{Preliminary comments on representations}

Above we have considered a~representation of the real Lie algebra $\overset{\longrightarrow}{\mathcal W}$ by real
derivations on a~space of complex-valued functions on $\Loop^1(\mathbb C\setminus\{0\})$.
This representation is real, in the sense that the set of real functions is stable, or equivalently that the action
commutes with complex conjugation of functions.

To be precise, f\/ix $\lambda\in \mathbb C$.
The Duren--Schif\/fer formulas imply that there is a~real representation of the real Lie algebra
$\overset{\longrightarrow}{\mathcal W}$ by real derivations on the spaces of complex-valued functions
\begin{gather*}
\mathbb C\big[u_1,\overline{u_1},u_2,\dots;\rho_0,\rho_0^{-1}\big]\rho_0^{\lambda},
\qquad
\mathbb C\big[l_1,\overline{l_1},l_2,\dots;\rho_{\infty},\rho_{\infty}^{-1}\big]\rho_{\infty}^{-\lambda},
\end{gather*}
and
\begin{gather*}
\mathbb C\big[u_1,\overline{u_1},u_2,\dots;l_1,\overline{l_1},l_2,\dots;\rho_0,\rho_0^{-1};\rho_{\infty},\rho_{\infty}^{-1}\big]a^{\lambda}.
\end{gather*}

Denote this real action of $\overset{\longrightarrow}{\mathcal W}$ by $\pi_0$.
By abstract nonsense there is an associated complex representation of $\mathcal W$ by complex derivations of the algebra
of complex-valued functions of self-avoiding loops, def\/ined~by
\begin{gather*}
\pi(L)=\frac12 \big(\pi_0\big(\overset{\rightarrow}{L}\big)-i\pi_0\big(\overset{\longrightarrow}{iL}\big)\big).
\end{gather*}
There is also a~representation
\begin{gather*}
\overline{\pi}\big(\overline{L}\big)=\frac12\big(\pi_0\big(\overset{\rightarrow}{L}\big)+i\pi_0\big(\overset{\longrightarrow}{iL}\big)\big).
\end{gather*}
This is a~complex representation of $\overline{\mathcal W}=\mathcal W^*$ by complex derivations.

In turn, in terms of the real embedding~\eqref{realembedding}
\begin{gather*}
\pi_0\big(\overset{\rightarrow}{L}\big)=\pi(L)+\overline{\pi}\big(\overline{L}\big).
\end{gather*}

The point of this translation is that the complex representations~$\pi$ and $\overline{\pi}$ are easier to analyze.
In fact (on proper domains) they can be expressed in terms of highest weight representations, and this allows us to
access well-known results from the theory of highest weight representations of the Virasoro algebra (at the moment the
central charge $c=0$, so that we are only considering the Witt algebra).

\subsection[Formulas for the representation~$\pi$]{Formulas for the representation~$\boldsymbol{\pi}$}

\begin{Proposition}\label{variation1}\quad
\begin{enumerate}\itemsep=0pt
\item[$(a)$] $(\pi(L_0)\phi_+)(z)=\phi_+(z)-\frac12 z\phi_+'(z)$.
In particular
\begin{gather*}
\pi(L_0)\rho_0=\frac12 \rho_0
\qquad
\text{and}
\qquad
\pi(L_0) u_k=-\frac12 ku_k,
\qquad
k\ge1.
\end{gather*}

\item[$(b)$] For $n>0$, $\pi(L_n)\phi_+=\phi_+^{n+1}$.
In particular
\begin{gather*}
\pi(L_n)\rho_0=0,
\qquad
\pi(L_n) u_k=0
\qquad
k<n,
\qquad
\pi(L_n) u_n=\rho_0^n
\end{gather*}
and in general
\begin{gather*}
\pi(L_n) u_{k}=\rho_0^n p_{k-n}^{(n+1)}(u_1,u_2,\dots),
\end{gather*}
where
\begin{gather*}
\big(1+u_1z+u_2z^2+\cdots\big)^{n+1}=\sum\limits_{l=0}^{\infty}p_l^{(n+1)}z^l.
\end{gather*}

\item[$(c)$] $\pi(L_{-1})\phi_+=1+u'(z)(-1+u_1 z)$.
In particular
\begin{gather*}
\pi(L_{-1})\rho_0=-u_1,
\qquad
\pi(L_{-1})(\rho_0u_1)=-3u_2+2u_1^2
\end{gather*}
and in general
\begin{gather*}
\pi(L_{-1})(\rho_0u_n)=-(n+2)u_{n+1}+(n+1)u_n u_1.
\end{gather*}
Hence
\begin{gather*}
\pi(L_{-1})u_n=\frac{n+2}{\rho_0}(u_1u_n-u_{n+1}).
\end{gather*}

\item[$(d)$] For $n>1$, $\pi(L_{-n})\rho_0=\frac12\rho_0^{-n+1}P_n(u_1,\dots,u_n)$, where
\begin{gather*}
P_n(u_1,\dots,u_n)=B_0(n)=\Res\left(\left(\frac{U'(t)}{U(t)}\right)^2t^{-n+1},t=0\right).
\end{gather*}
If $\deg(u_j)=j$, then $P_n$ is a~homogeneous polynomial of degree~$n$.

\item[$(e)$] For $k\ge1$
\begin{gather*}
\pi(L_{-n})u_k=\rho_0^{-n}\left(\frac{k}{2}u_kB_0(n)+ \sum\limits_{m=1}^{k}(k+1-m)u_{k-m}B_{-m}(n)\right).
\end{gather*}
Equivalently
\begin{gather*}
\pi(L_{-n})(\rho_0 u_k)=\rho_0^{-n+1}\sum\limits_{m=0}^{k}(k+1-m)u_{k-m}B_{-m}(n)-\rho_0^{-n+1}\frac{k+1}{2}u_kB_0(n),
\end{gather*}
where
\begin{gather*}
B_m(n)=\Res\left(\left(\frac{U'(t)}{U(t)}\right)^2 U(t)^m t^{-n+1},t=0\right)=\Res\left(\frac{u(z)^{-n+1}}{u'(z)}z^{m-2},z=0\right).
\end{gather*}
\end{enumerate}

\end{Proposition}

Using Lemma~\ref{lemma1} below, this can be restated in the following way.

\begin{Proposition}\label{restate1}\quad
\begin{enumerate}\itemsep=0pt
\item[$(a)$] For $n\in\mathbb Z$
\begin{gather*}
\pi(L_{n})\rho_0=\frac12\rho_0 \Res\left(\frac{\phi_+^{n+1}(z)}{z^2\phi_+'(z)},z=0\right).
\end{gather*}

\item[$(b)$] For $k\ge1$
\begin{gather*}
\pi(L_{n})u_k=\frac{k}{2}u_k\widetilde{B}_0(n)+ \sum\limits_{m=1}^{k}(k+1-m)u_{k-m}\widetilde{B}_{-m}(n).
\end{gather*}
Equivalently
\begin{gather*}
L_{n}(\rho_0u_k)=\rho_0^{n+1}\sum\limits_{m=0}^{k}(k+1-m)u_{k-m}B_{-m}(n)-\rho_0^{n+1}\frac{k+1}{2}u_kB_0(n),
\end{gather*}
where
\begin{gather*}
\widetilde{B}_m(n)=\Res\left(\frac{\phi_+(z)^{n+1}}{\phi_+'(z)} z^{m-2},z=0\right).
\end{gather*}
\end{enumerate}
\end{Proposition}

\begin{Remark}
This second statement seems cleaner than the f\/irst.
However, as we will see when we introduce the energy-momentum tensor, the f\/irst statement has the advantage of being
stated in terms of the inverse of $\phi_+$.
\end{Remark}

To avoid cumbersome notation, we will often identify $L_n$ with its corresponding operator, $\pi(L_n)$.
Suppose that we write $u_0=1$ and $a_k=\rho_0u_k$, so that
\begin{gather*}
\phi_+(z)=\sum\limits_{k=0}^{\infty}a_kz^{k+1}.
\end{gather*}

If $n>0$ and $k\ge1$, then according to (e)
\begin{gather*}
L_{-n}(a_k)=\sum\limits_{m=0}^{k}(k+1-m)a_{k-m} \Res\left(\left(\frac{\big(\phi_+^{-1}\big)'(t)}{\phi_+^{-1}(t)}\right)^2
\phi_+^{-1}(t)^{-m} t^{-n+1},t=0\right)
\\
\phantom{L_{-n}(a_k)=}{}
-\frac{k+1}{2}a_{k} \Res\left(\left(\frac{\big(\phi_+^{-1}\big)'(t)}{\phi_+^{-1}(t)}\right)^2 t^{-n+1},t=0\right)
\end{gather*}
and
\begin{gather*}
L_{-n}(\phi_+)=\sum\limits_{k=0}^{\infty}\left(\sum\limits_{m=0}^{k}(k+1-m)a_{k-m}
\Res\left(\left(\frac{\big(\phi_+^{-1}\big)'(t)}{\phi_+^{-1}(t)}\right)^2 \phi_+^{-1}(t)^{-m} t^{-n+1},t=0\right)
\right.
\\
\left.\phantom{L_{-n}(\phi_+)=}{}
-\frac12(k+1)a_{k} \Res\left(\left(\frac{\big(\phi_+^{-1}\big)'(t)}{\phi_+^{-1}(t)}\right)^2 t^{-n+1},t=0\right)\right)z^{k+1}.
\end{gather*}

\begin{Lemma}\label{lemma1}
\begin{gather*}
\Res\left(\left(\frac{\big(\phi_+^{-1}\big)'(t)}{\phi_+^{-1}(t)}\right)^2 \phi_+^{-1}(t)^{-m}t^{-n+1},t=0\right)
=\Res\left(\frac{\phi_+(z)^{-n+1}}{\phi_+'(z)z^{2+m}},z=0\right).
\end{gather*}
\end{Lemma}

\begin{proof}
Fix a~small circle~$C$ surrounding $0$ in the~$t$ plane.
Then
\begin{gather*}
\int_C\left(\frac{\big(\phi_+^{-1}\big)'(t)}{\phi_+^{-1}(t)}\right)^2 \phi_+^{-1}(t)^{-m}
t^{-n+1}dt=\int_{\phi_+^{-1}(C)} \big(\phi_+^{-1}\big)'(\phi_+(z))^2\frac{\phi_+(z)^{-n+1}}{z^{m+2}}d\phi_+(z)
\\
\hphantom{\int_C\left(\frac{\big(\phi_+^{-1}\big)'(t)}{\phi_+^{-1}(t)}\right)^2 \phi_+^{-1}(t)^{-m}t^{-n+1}dt}{}
=2\pi i \Res\left(\frac{\phi_+(z)^{-n+1}}{\phi_+'(z)z^{2+m}},z=0\right).\tag*{\qed}
\end{gather*}
\renewcommand{\qed}{}
\end{proof}

This can be restated more cleanly in the following way.

\begin{Lemma}
\begin{gather*}
\rho_0^{-n}B_m(n)=\Res\left(\frac{\phi_+^{-n+1}(z)}{\phi_+'(z)}z^{m-2},z=0\right).
\end{gather*}
\end{Lemma}

Using the lemma we can write
\begin{gather*}
L_{-n}(\phi_+)=\sum\limits_{k=0}^{\infty}\left(\sum\limits_{m=0}^{k}(k+1-m)a_{k-m}
\Res\left(\frac{\phi_+(s)^{-n+1}}{\phi_+'(s)s^{2+m}},s=0\right)
\right.
\\
\left.\phantom{L_{-n}(\phi_+)=}{}
-\frac12(k+1)a_{k} \Res\left(\frac{\phi_+(s)^{-n+1}}{\phi_+'(s)s^2},s=0\right)\right)z^{k+1}
\\
\phantom{L_{-n}(\phi_+)}{}
=\phi_+'(z)\left(\sum\limits_{m=0}^{\infty}\left(\frac{\phi_+^{-n+1}}{\phi_+'}\right)_{m+1}z^{m+1}
-z\Res\left(\frac{\phi_+(s)^{-n+1}}{\phi_+'(s)s^2},s=0\right)\right)
\\
\phantom{L_{-n}(\phi_+)}{}
=\phi_+'(z)\left(\left(\frac{\phi_+(z)^{-n+1}}{\phi_+'(z)}\right)_{++} -z\Res\left(\frac{\phi_+(s)^{-n+1}}{\phi_+'(s)s^2},s=0\right)\right).
\end{gather*}

The pleasant surprise is that this expression leads to a~formula which is valid for all~$n$.

\begin{Theorem}
\label{uniformformula1}
For any $n\in\mathbb Z$
\begin{gather*}
L_{n}(\phi_+)(z)
=\phi_+'(z)\left(\frac{\phi_+(z)^{n+1}}{\phi_+'(z)}\right)_{++} -\frac12 z\phi_+'(z)\Res\left(\frac{\phi_+(s)^{n+1}}{\phi_+'(s)s^2},s=0\right)
\end{gather*}
and
\begin{gather*}
L_{n}(u)(z)=\rho_0^{n}\left(u'(z)\left(\frac{u(z)^{n+1}}{u'(z)}\right)_{++} -\frac12(zu'(z)+u(z))\Res\left(\frac{u(s)^{n+1}}{u'(s)s^2},s=0\right)\right).
\end{gather*}

\end{Theorem}

\begin{proof}
We just need to check that this formula agrees with our previous calculations when $n\ge0$.
This is straightforward.
\end{proof}

\subsection[Formulas for $\overline{\pi}$]{Formulas for $\boldsymbol{\overline{\pi}}$}

\begin{Proposition}\label{variation2}\quad
\begin{enumerate}\itemsep=0pt
\item[$(a)$] $\overline{\pi}(\overline{L}_0)\phi_+=\frac12 z\phi_+'(z)$.
In particular
\begin{gather*}
\overline{\pi}(\overline{L}_0)\rho_0=\frac12 \rho_0
\qquad
\text{and}
\qquad
\overline{\pi}(\overline{L}_0) u_k=\frac{k-1}{2} u_k,
\qquad
k\ge1.
\end{gather*}

\item[$(b)$] For $n>0$, $\overline{\pi}(\overline{L}_n)\phi_+=0$.

\item[$(c)$] $\overline{\pi}(\overline{L}_{-1})\phi_+=u'(z)(-u_1^* z+z^2)$.
In particular
\begin{gather*}
\overline{\pi}(\overline{L}_{-1})\rho_0=-u_1^*,
\qquad
\overline{\pi}(\overline{L}_{-1})u_1=\rho_0^{-1}(1-u_1u_1^*).
\end{gather*}
In general
\begin{gather*}
\overline{\pi}(\overline{L}_{-1})u_n=\rho_0^{-1}n(u_{n-1}-u_1^*u_n).
\end{gather*}

\item[$(d)$] For $n>1$, $\overline{\pi}(\overline{L}_{-n})\rho_0=\frac12\rho_0^{-n+1}P_n(u_1,\dots,u_n)^*$.

\item[$(e)$] For $k\ge1$
\begin{gather*}
\overline{L}_{-n}u_k=\rho_0^{-n}\left(\frac{k}{2}u_kB_0(n)^*+ \sum\limits_{m=1}^{k}(k+1-m)u_{k-m}B_{m}(n)^*\right).
\end{gather*}
Equivalently
\begin{gather*}
\overline{L}_{-n}(\rho_0u_k)=\rho_0^{-n+1}\sum\limits_{m=0}^{k}(k+1-m)u_{k-m}B_{m}(n)^*-\rho_0^{-n+1}\frac{k+1}{2}u_kB_0(n)^*.
\end{gather*}
\end{enumerate}
\end{Proposition}

Now we want to add things up as in the preceding section.
As before we write $\phi_+(z)=\sum a_kz^{k+1}$, where $a_k=\rho_0u_k$ and it is understood that $u_0=1$.
By part (e)
\begin{gather*}
\overline{L}_{-n}(a_k)=\rho_0^{-n}\sum\limits_{m=0}^{k}(k+1-m)a_{k-m}B_{m}(n)^*-\rho_0^{-n}\frac{k+1}{2}a_kB_0(n)^*.
\end{gather*}
By the change of variable lemma of the preceding subsection
\begin{gather*}
\rho_0^{-n}B_m(n)=\Res\left(\frac{\phi_+(z)^{-n+1}}{\phi_+'(z)}z^{m-2},z=0\right).
\end{gather*}
Therefore
\begin{gather*}
\overline{L}_{-n}\phi_+(z)=\phi_+'(z)\sum\limits_{m=0}^{\infty}\left(\left(\frac{\phi_+(z)^{-n+1}}{\phi_+'(z)}\right)_{-m+1}\right)^*z^{m+1}
-\frac12z\phi_+'(z)\left(\left(\frac{\phi_+(z)^{-n+1}}{\phi_+'(z)}\right)_1\right)^*,
\end{gather*}
where the notation $(\cdots)_k$ denotes the $k$th Fourier coef\/f\/icient.
This equals
\begin{gather*}
\phi_+'(z)\sum\limits_{m=0}^{\infty}\left[\left(\frac{\phi_+(z)^{-n+1}}{\phi_+'(z)}\right)^*z^2\right]_{m+1} z^{m+1}
-\frac12z\phi_+'(z)\left(\left(\frac{\phi_+(z)^{-n+1}}{\phi_+'(z)}\right)_1\right)^*
\\
\qquad{}
=\left[\left(z^{-2}\frac{\phi_+^{-n+1}(z)}{\phi_+'(z)}\right)^*\right]_{++}
-\frac12z\phi_+'(z)\left(\left(\frac{\phi_+(z)^{-n+1}}{\phi_+'(z)}\right)_1\right)^*.
\end{gather*}

As in the preceding subsection, we obtain the following uniform formula.

\begin{Theorem}
\label{uniformformula2}
For any $n\in \mathbb Z$
\begin{gather*}
\overline{L}_{n}\phi_+(z)=\phi_+'(z)\left(\left[\left(z^{-2}\frac{\phi_+^{n+1}(z)}{\phi_+'(z)}\right)^*\right]_{++}
-\frac12z\left((\frac{\phi_+(z)^{n+1}}{\phi_+'(z)})_1\right)^*\right).
\end{gather*}
\end{Theorem}

\begin{proof}
We just need to check that this formula agrees with the formulas in Proposition~\ref{variation2}.
This is again straightforward.
\end{proof}

\subsection[Formulas for $\pi_0$, revisited]{Formulas for $\boldsymbol{\pi_0}$, revisited}

We can use Theorems~\ref{uniformformula1} and~\ref{uniformformula2} to recast the Duren--Schif\/fer variational formulas
in the following form.

\begin{Corollary}
For all $n\in \mathbb Z$, $\overset{\rightarrow}{L}_n\phi_+$ equals
\begin{gather*}
\phi_+'(z) \left[\frac{\phi_+(z)^{n+1}}{\phi_+'(z)}+\left(z^{-2}\frac{\phi_+(z)^{n+1}}{\phi_+'(z)}\right)^*\right]_{++}\\
\qquad \quad{} -\frac12 z\left(\Res\left(\frac{\phi_+(s)^{n+1}}{\phi_+'(s)s^2},s=0\right)
+\Res\left(\frac{\phi_+(s)^{n+1}}{\phi_+'(s)s^2},s=0\right)^*\right)\\
\qquad{}
=\phi_+'(z)\left(\frac12(c_1+\overline{c_1})z+\sum\limits_{k>1}(c_k+\overline{c_{2-k}})z^k\right),
\end{gather*}
where
\begin{gather*}
\frac{\phi_+(z)^{n+1}}{\phi_+'(z)}=\sum\limits_{k=n+1}^{+\infty} c_kz^k.
\end{gather*}
\end{Corollary}

\begin{proof}
By def\/inition
\begin{gather*}
\overset{\longrightarrow}{L}_n\phi_+=L_n\phi_+ +\overline{L}_n\phi_+.
\end{gather*}
Theorems~\ref{uniformformula1} and~\ref{uniformformula2} imply that this equals
\begin{gather*}
\phi_+'(z)\left[\frac{\phi_+(z)^{n+1}}{\phi_+'(z)}\right]_{++} -\frac12 z\phi_+'(z)\Res\left(\frac{\phi_+(s)^{n+1}}{\phi_+'(s)s^2},s=0\right)
\\
\qquad{}
+\phi_+'(z)\left[\left(z^{-2}\frac{\phi_+(z)^{n+1}}{\phi_+'(z)}\right)^*\right]_{++}
-\frac12z\phi_+'(z)\left(\left(\frac{\phi_+(z)^{n+1}}{\phi_+'(z)}\right)_1\right)^*.\tag*{\qed}
\end{gather*}
\renewcommand{\qed}{}
\end{proof}

It is obviously desirable to f\/ind a~direct proof of these formulas which ref\/lects their structure.

\subsection[Calculations with $\phi_-$]{Calculations with $\boldsymbol{\phi_-}$}

On the one hand, in the standard~$w$ coordinate at $\infty \in \mathbb P^1$,
\begin{gather*}
\frac{1}{\phi_-(\frac{1}{w})}=\frac1{\rho_{\infty}}w\left(1+\sum\limits_{n=1}^{\infty}l_nw_n\right).
\end{gather*}
The $l_n$ coordinates for $\phi_-$ are analogous to the $u_n$ coordinates for $\phi_+$, and variational formulas for
$\phi_-$ essentially arise from substituting $l_j$'s for $u_j$'s in our earlier formulas.
On the other hand, in the standard~$z$ coordinate,
\begin{gather*}
\phi_-(z)=\rho_{\infty}L(z)=\rho_{\infty}z\left(1+\sum\limits_{m=1}^{\infty}b_mz^{-m}\right)
\end{gather*}
and it is occasionally useful to employ the $b_m$ coordinates.
The relation between the two sets of coordinates is standard.

\begin{Lemma}
\begin{gather*}
\mathbb C[l_1,l_2,\dots]=\mathbb C[b_1,b_2,\dots].
\end{gather*}
In fact for each~$M$
\begin{gather*}
\mathbb C[l_1,l_2,\dots,l_M]=\mathbb C[b_1,b_2,\dots,b_M].
\end{gather*}
\end{Lemma}

\begin{proof}
\begin{gather*}
w\left(1+\sum\limits_{n=1}^{\infty}l_nw^n\right)=\frac{1}{\frac{1}{w}\left(1+\sum\limits_{m=1}^{\infty}b_mw^{m}\right)}
\end{gather*}
or
\begin{gather*}
 1+\sum\limits_{n=1}^{\infty}l_nw^n =\frac{1}{1+\sum\limits_{m=1}^{\infty}b_mw^{m}}
\end{gather*}
implies
\begin{gather*}
l_1=-b_1,
\qquad
l_2=-b_2+b_1^2,
\qquad
\dots.\tag*{\qed}
\end{gather*}
\renewcommand{\qed}{}
\end{proof}

The $\phi_-$ analog of Theorems~\ref{uniformformula1} and~\ref{uniformformula2} is the following theorem.
In the statement, for a~Laurent expansion convergent in an annulus $R<\vert z\vert<\infty$, we use the notation
$\Res(\sum g_mz^m,z=\infty)=-g_{-1}$ (This is actually the residue of the dif\/ferential $g(z)dz$ at $z=\infty$ in the Riemann sphere).

\begin{Theorem}
Let $n\in \mathbb{Z}$.
\begin{enumerate}\itemsep=0pt

\item[$(a)$]
\begin{gather*}
L_n (\phi_-(z) )=-z^2\phi_-'(z)
\left[\left(\frac{\phi_-(z)^{n+1}}{z^2\phi_-'(z)}\right)_{-}+\frac{z^{-1}}{2}\Res \left(\frac{\phi_-(t)^{n+1}}{t^2\phi_-'(t)},
t=\infty\right)\right]
\end{gather*}
and
\begin{gather*}
L_n (L(z) )=-\frac{\rho_{\infty}^n}{2} \Res \left(\frac{L(t)^{n+1}}{t^2L'(t)}, t=\infty \right)L(z)\\
\hphantom{L_n (L(z) )=}{}
-\rho_{\infty}^nz^2L'(z)\left[\left(\frac{L(z)^{n+1}}{z^2L'(z)}\right)_{-}+\frac{z^{-1}}{2}\Res \left(\frac{L(t)^{n+1}}{t^2L'(t)},
 t=\infty\right)\right].
\end{gather*}

\item[$(b)$]
\begin{gather*}
\overline{L}_n (\phi_-(z) )=-z^2\phi_-'(z)
\left[\frac{z^{-1}}{2}
\Res\left(\frac{\phi_-(t)^{n+1}}{t^2\phi_-'(t)}, \infty\right)^{\ast} +\left(\left(\frac{\phi_-(z)^{n+1}}{\phi_-'(z)}\right)^{\ast}\right)_-\right]
\end{gather*}
and
\begin{gather*}
\overline{L}_n (L(z) )=-\frac{\rho_{\infty}^n}{2}\big(L(z)+zL'(z)\big)\Res \left(\frac{L(t)^{n+1}}{t^2L'(t)},
t=\infty\right)^{\ast}\\
\hphantom{\overline{L}_n (L(z) )=}{}
-\rho_{\infty}^nz^2L'(z) \left(\left(\frac{L(z)^{n+1}}{L'(z)} \right)^{\ast} \right)_-.
\end{gather*}
\end{enumerate}
\end{Theorem}

\subsection{Representation-theoretic consequences}
\label{representations}

The formulas of the preceding section imply that~$\pi$ is a~complex representation of the Witt algebra $\mathcal W$~by
derivations of the algebra $\Omega^0(\rho_0)\otimes \mathbb C[u_1,u_2,\dots]$, where $\Omega^0(\rho_0)$ denotes any algebra
of smooth functions of $\rho_0$.

Consider the action of $\mathcal W$ on the vector space
\begin{gather*}
\mathbb C\big[\rho_0^{\lambda},\rho_0,\rho_0^{-1},u_1,u_2,\dots\big],
\end{gather*}
where~$\lambda$ is a~f\/ixed complex number.
For $n>0$ the operators $L_n$ kill $\rho_0^{\lambda}$, and the spectrum of $L_0$ on the $\mathcal W$-module generated~by
$\rho_0^{\lambda}$ is $\{\lambda/2+n: n=0,1,\dots\}$.
We will refer to this as a~lowest weight module (admittedly there are conf\/licting conventions).
The following proposition follows from well-known facts about such representations (see~\cite{Kac}).

\begin{Proposition}
For any $\lambda\in \mathbb C$,
\begin{enumerate}\itemsep=0pt
\item[$(a)$] The representation generated by the~$\pi$ action of $\mathcal W$ on $\rho_0^{\lambda}$ is a~realization of the
unique irreducible lowest weight representation of the Virasoro algebra with central charge $c=0$ and
$h=\frac12\lambda$.
If $\lambda\ne -\frac{m^2-1}{12}$, then
\begin{gather*}
\pi(\mathcal U(\mathcal W))\rho_0^{\lambda}=\bigoplus_{n=0}^{\infty}\rho_0^{\lambda-n}\mathbb C[u_1,u_2,\dots]^{(n)},
\end{gather*}
where $u_j$ has degree~$j$.
Otherwise there is a~proper containment.

\item[$(b)$] Similarly, the representation generated by the~$\pi$ action of $\mathcal W$ on $\rho_{\infty}^{-\lambda}$ is
a~realization of the highest weight representation of the Virasoro algebra with central charge $c=0$ and
$h=-\frac{\lambda}{2}$.
If $\lambda\ne -\frac{m^2-1}{12}$, then
\begin{gather*}
\pi(\mathcal U(\mathcal W))\rho_{\infty}^{-\lambda}=\bigoplus_{n=0}^{\infty}\rho_{\infty}^{-\lambda-n}\mathbb
C[l_1,l_2,\dots]^{(n)},
\end{gather*}
where $l_j$ has degree~$j$.
Otherwise there is a~proper containment.
\end{enumerate}
\end{Proposition}

\begin{Remark}
The realization of the lowest weight representation in part (a) is related in a~relatively simple way to the
realization, using geometric quantization techniques, due to Kirillov and Yuriev in~\cite{KY}.
In~\cite{KY} $\mathcal W$ acts on a~space of sections of a~line bundle (parameterized by $c=0$ and $h=\lambda/2$) over
(a somewhat imprecisely def\/ined) space of Schlicht functions $u\in \mathcal S$ (normalized univalent functions on the
disk, viewed as a~homogeneous space for $\Diff(S^1)$).
In coordinates (by trivializing the line bundle) this vector space is identif\/ied with $\mathbb C[u_1,u_2,\dots]$,
polynomials in the coef\/f\/icients of the univalent function~$u$, and the formulas for the action appear in~(8)
of~\cite{KY} (with $c=0$, and one takes the negative of the operators, because we consider the opposite of the bracket
in~\cite{KY}).
The intertwining operator from Kirillov and Yuriev's realization to our realization in (a) is given by the map
\begin{gather*}
\mathbb{C}[u_1,u_2,\dots ]\to \mathbb C\big[\rho_0^{\lambda},\rho_0,\rho_0^{-1},u_1,u_2,\dots\big]: \ P(u_1,u_2,\dots) \mapsto
P\big(U_1/\rho_0,U_2/\rho_0^2,\dots \big)\rho_0^{\lambda},
\end{gather*}
where $U=t(1+\sum\limits_{n>0}U_nt^n)$ is the inverse to the univalent function $u=z\big(1+\sum\limits_{n>0}u_nz^n\big)$.
An advantage of our realization is that the operators are derivations of an algebra, which makes them more amenable to
calculations.
This will appear in the f\/irst author's dissertation.
\end{Remark}

\subsection{Stress-energy formulation}

Consider the standard holomorphic coordinate $z=x+iy$.
In real coordinates the symmetric stress tensor has the form
\begin{gather*}
\mathcal T=\left(
\begin{matrix}
dx &dy
\end{matrix}
\right)\left(
\begin{matrix}
T_{11}&T_{12}
\\
T_{21}&T_{22}
\end{matrix}
\right)\left(
\begin{matrix}
dx
\\
dy
\end{matrix}
\right),
\end{gather*}
where $T_{12}=T_{21}$.
In complex coordinates
\begin{gather*}
T=\left(
\begin{matrix}
dz &d\overline z
\end{matrix}
\right)\left(
\begin{matrix}
T_{11}-T_{22}-2iT_{12} &T_{11}+T_{22}+i(T_{12}-T_{21})
\\
T_{11}-T_{22}+i(T_{21}-T_{12})&T_{11}-T_{22}+2iT_{12}
\end{matrix}
\right)\left(
\begin{matrix}
dz
\\
d\overline{z}
\end{matrix}
\right)\\
\hphantom{T}{}
=\left(
\begin{matrix}
dz &d\overline z
\end{matrix}
\right)\left(
\begin{matrix}
T_{11}-T_{22}-2iT_{12} &T_{11}+T_{22}
\\
T_{11}+T_{22}&T_{11}+T_{22}+2iT_{12}
\end{matrix}
\right)\left(
\begin{matrix}
dz
\\
d\overline{z}
\end{matrix}
\right).
\end{gather*}
Conformal invariance is implied by the trace condition
\begin{gather*}
\operatorname{tr}(T)=T_{11}+T_{22}=0
\end{gather*}
(see~\cite[p.~101 and p.~103]{DMS}).
In complex coordinates this implies that~$T$ is diagonal.

In a~conformal f\/ield theory with central charge $c=0$
\begin{gather*}
T(z):=(T_{11}-T_{22}-2iT_{12})dz^2=\sum\limits_{n=-\infty}^{\infty}L_nz^{-n}\left(\frac{dz}{z}\right)^2
\end{gather*}
is a~holomorphic quadratic dif\/ferential (see~\cite[p.~155]{DMS};
note: for $c\ne0$, the stress energy ``tensor'' is actually a~holomorphic projective connection; see~\cite[p.~136]{DMS} or~\cite[p.~532]{Segal}).

We are seeking a~completely natural formulation for the action of the Witt algebra

\begin{Proposition}
\begin{gather*}
T(t)\rho_0=\frac{\rho_0}{2}\big(\partial \log \big(\phi_+^{-1}(t)\big)\big)^2,
\qquad
T(t)\rho_{\infty}=-\frac{\rho_{\infty}}{2}\big(\partial \log \big(\phi_-^{-1}(t)\big)\big)^2.
\end{gather*}
\end{Proposition}

\begin{proof}
By def\/inition
\begin{gather*}
T(t)\rho_0=\sum\limits_{n=-\infty}^{\infty}L_{-n}(\rho_0)t^n\left(\frac{dt}{t}\right)^2.
\end{gather*}
By part (a) of Proposition~\ref{restate1}, this equals
\begin{gather*}
\frac12\rho_0\sum\limits_{n=0}^{\infty} \Res\left(\left(\frac{\big(\phi_+^{-1}\big)'(t)}{\phi_+^{-1}(t)}\right)^2
t^{-n+1},t=0\right)t^n\left(\frac{dt}{t}\right)^2=\frac{\rho_0}{2}\big(\partial \log \big(\phi_+^{-1}(t)\big)\big)^2 .
\end{gather*}

This proves the f\/irst statement.
The proof of the second statement is similar.
\end{proof}

\begin{Corollary}
In the sense of hyperfunctions
\begin{gather*}
T(t)a^{\lambda}=\frac{\lambda}{2} \left(\big(\partial \log\big(\phi_+^{-1}(t)\big)\big)^2+\big(\partial
\log\big(\phi_-^{-1}(t)\big)\big)^2\right)a^{\lambda}.
\end{gather*}
\end{Corollary}

\begin{proof}
From a~formal power series point of view, this follows immediately from the proposition.
From the point of view of analysis, this equality has to be interpreted in a~hyperfunction sense, because the f\/irst term
is holomorphic in $U_+$ and the second term is holomorphic in~$U_-$.
\end{proof}

\section{Inf\/initesimal invariance}
\label{infinitesimalinvariance}

Suppose that $\gamma\in \Loop^1(\mathbb C \setminus\{0\})$.
In terms of the standard coordinate~$z$,
\begin{gather*}
\phi_+(z)=\rho_0(\gamma)u(z),
\qquad
u(z)=z\bigg(1+\sum\limits_{n\ge1}u_nz^{n}\bigg).
\end{gather*}
In terms of the coordinate $w=\frac1{z}$,
\begin{gather*}
\frac1{\phi_-\big(\frac1{w}\big)}=\frac1{\rho_{\infty}}w\bigg(1+\sum\limits_{n\ge1}l_nw^{n}\bigg).
\end{gather*}

The variational formulas of the preceding section imply that the vector space of functions of the form
\begin{gather*}
p(u_1,\dots,u_n,\overline{u_1},\dots,\overline{u_n})f(\rho_0),
\end{gather*}
where~$p$ is a~polynomial of any number of variables, and~$f$ has compact support in $\mathbb R^+$, is stable with
respect to the action of the Witt algebra (this applies both to the real action and the complexif\/ied actions).
Since the Witt algebra is stable with respect to $z\mapsto w=\frac1{z}$, the vector space of functions of the form
\begin{gather*}
p\big(l_1,\dots,l_n,\overline{l_1},\dots,\overline{l_n}\big)f(\rho_{\infty}),
\end{gather*}
where~$p$ is a~polynomial of any number of variables, and~$f$ has compact support in $\mathbb R^+$, is also stable with
respect to the action of the Witt algebra.
Consequently the vector space of ``test functions'' spanned by functions of the form
\begin{gather}
\label{testfunction}
F=p\big(u_1,\dots,\overline{u_1},\dots,l_1,\dots,\overline{l_1},\dots \big)f(\rho_0,\rho_{\infty}),
\end{gather}
where~$p$ is a~polynomial and~$f$ has compact support in $\mathbb R^+\times \mathbb R^+$, is stable with respect to the
Witt algebra (for the real or complexif\/ied actions).
In reference to~$F$, since $u_n$ and $l_n$ are bounded (by constants depending only on~$n$),~$p$ is bounded.
The compact support condition on~$f$ implies that~$F$ is supported on $\Loop^1$ of a~f\/ixed f\/inite type annulus.
Since $\mu_0$ has f\/inite measure on loops in a~f\/inite type annulus,~$F$ is integrable.

\begin{Proposition}
\label{intbyparts1}
The measure $\mu_0$ is infinitesimally conformally invariant, in the sense that for any $L\in\mathcal W\times
\overline{\mathcal W}$
\begin{gather*}
\int L(F)d\mu_0(\gamma)=0
\end{gather*}
for any test function~$F$ as in~\eqref{testfunction}.

\end{Proposition}

\begin{proof}
It suf\/f\/ices to prove the proposition for $\overset{\rightarrow}{L}\in \overset{\longrightarrow}{\mathcal W}$.

By Koebe's theorem, a~test function~$F$ as in~\eqref{testfunction} is supported on $\Loop^1(\{\delta<\vert
z\vert<\delta^{-1}\})$ for some~$\delta$.
Let $A_0$ denote a~f\/inite type annulus containing $\{\delta\le\vert z\vert\le\delta^{-1}\}$.
For some positive~$t_0$, for all $\vert t\vert<t_0$, the f\/low $\exp (t\vec{L})$ is def\/ined on $A_0$, and
$\{\delta\le\vert z\vert\le\delta^{-1}\}$ will be contained in $\cap_{t<t_0}A_t$, where $A_t:=\exp (t\vec{L})A_0$.
By local conformal invariance
\begin{gather*}
\int F d\mu_0=\int F d\mu_{A_0}=\int \big(e^{t\vec{L}}\big)_*F d\mu_{A_t}=\int \big(e^{t\vec{L}}\big)_*F d\mu_0.
\end{gather*}

To complete the proof we need to justify taking the derivative with respect to~$t$ at $t=0$ under the last integral.
The derivative $\vec{L}F$ is another test function, necessarily bounded.
The translates
\begin{gather}
\label{translates}
\big(e^{t\vec{L}}\big)_*\big(\vec{L}F\big),
\qquad
\vert t\vert <t_0
\end{gather}
are also uniformly bounded by the same constant.
Moreover the translates~\eqref{translates} are all supported on $\Loop^1$ of some f\/inite type annulus, for which
the $\mu_0$ measure is f\/inite.
Thus a~multiple of the characteristic function of $\Loop^1$ for this f\/ixed f\/inite type annulus is integrable and
dominates all of the translates~\eqref{translates}.
Hence by dominated convergence we can dif\/ferentiate under the integral sign.
\end{proof}

Kontsevich and Suhov have conjectured that there is a~converse of this result which holds generally for their
conjectural family of measures $\mu_c$ deforming $\mu_0$ (see Section~2.5.2 of~\cite{KS}).

For the purposes of this paper, we need to be able to apply integration by parts to functions which involve the bounded
function $a^{\lambda}$ ($\lambda>0$), rather than a~function having compact support in $\rho_0$, $\rho_{\infty}$.
One complication is that for $\vec{L}\in \vec{\mathcal W}$,
\begin{gather*}
\vec{L}a^{\lambda}=\lambda a^{\lambda-1}\vec{L}(a)
\end{gather*}
is not necessarily bounded.

\begin{Lemma}
\label{intbyparts2}
Suppose that
\begin{gather*}
F=p(u,\overline{u},l,\overline{l})f(\rho_{\infty})a^{\lambda},
\end{gather*}
where~$p$ is a~polynomial and~$f$ has compact support in $\mathbb R^+$.
Then for any $L\in\mathcal W\times \overline{\mathcal W}$, for $\Reop(\lambda)$ sufficiently large,
\begin{gather*}
\int L(F)d\mu_0(\gamma)=0.
\end{gather*}
The same conclusion applies if we replace $f(\rho_{\infty})$ by $f(\rho_{0})$.
\end{Lemma}

\begin{proof}
Fix a~smooth positive function $g(\rho_0)$ having compact support for $\rho_0\in\mathbb R^+$ and identically~$1$ in
a~neighborhood of $\rho=1$.
By Proposition~\ref{intbyparts1}, for each $\delta>0$,
\begin{gather*}
\int L(g(\delta\rho_0)F)d\mu_0(\gamma)=0
\end{gather*}
or
\begin{gather*}
\delta\int g'(\delta\rho_0)L(\rho_0)Fd\mu_0+\int g(\delta\rho_0)L(F)d\mu_0(\gamma)=0.
\end{gather*}
Since~$g$ is f\/ixed, the f\/irst term goes to zero as $\delta\to0$.
We can apply dominated convergence to the second term, for suf\/f\/iciently large~$\lambda$ (so that the part of the
integrand not involving~$g$ is bounded, and hence the integral is well-def\/ined).
This implies the Lemma.
\end{proof}

\begin{Proposition}
Suppose that
$F=a^{\lambda}p(u,\overline{u},l,\overline{l})$,
where~$p$ is a~polynomial.
\begin{enumerate}\itemsep=0pt
\item[$(a)$] If $L=L_n$ {or} $\overline{L}_n$ with $n\le0$, then for sufficiently large $\Reop(\lambda)$
\begin{gather*}
\int L(F)\vert_{\rho_{\infty}=1}d\nu_0=0.
\end{gather*}

\item[$(b)$] If $L=L_n$ {or} $\overline{L}_n$ with $n\ge0$, then for sufficiently large $\Reop(\lambda)$
\begin{gather*}
\int L(F)\vert_{\rho_0=1}d\nu_0=0.
\end{gather*}
\end{enumerate}
\end{Proposition}

\begin{proof}
Suppose that $L=L_n$ or $\overline{L}_n$ with $n<0$.
Fix a~smooth family of functions $g_{\delta}(\rho_{\infty})$ which converges to the~$\delta$ function at
$\rho_{\infty}=1$.
Using $L_n(\rho_{\infty})=0$ and Lemma~\ref{intbyparts2},
\begin{gather*}
\int L(Fg_{\delta}(\rho_{\infty}))d\mu_0 =\int L(F)g_{\delta}(\rho_{\infty})d\mu_0=0.
\end{gather*}
Since $L(F)$ is bounded for suf\/f\/iciently large $\Reop(\lambda)$, the left hand side of the last equality converges
to $\int L(F)d\nu_0$ as $\delta \to0$.
This implies part~(a).

If $L=L_n$ or $\overline{L}_n$ with $n>0$, the same argument applies with $g_{\delta}(\rho_0)$ in place of $g_{\delta}(\rho_{\infty})$.

If $L=\overset{\rightarrow}{L}_0$, then $L(F)=0$.
We have previously observed that if $L=\overset{\longrightarrow}{iL}_0$, then~$L$ exponentiates to rotational symmetry
of $\mathbb C\setminus\{0\}$, and this corresponds to invariance of $\nu_0$ with respect to the conjugation action of
rotations on homeomorphisms.
\end{proof}

In the sections below, we will repeatedly apply a~variation of the preceding proof in the following way.
Suppose that $n>0$ and $L=L_n$ or $L=\overline{L}_n$.
Then as in the proof
\begin{gather*}
\int L\big(\rho^{-n}Fg_{\delta}(\rho_{0})\big)d\mu_0 =\int L(F)\rho_0^{-n}g_{\delta}(\rho_{0})d\mu_0=0.
\end{gather*}
We can take the limit as $\delta\to0$, because the support of $g_{\delta}$ remains bounded, and $\rho_0^{-n}$ will be
bounded in this support region.
This implies
\begin{gather*}
\int L(F)\rho_0^{-n}d\nu_0=0,
\end{gather*}
which can be written heuristically as
\begin{gather*}
\int L(F)\rho_0^{-n}\delta_1(\rho_{0})d\mu_0=0,
\end{gather*}
where $\delta_1$ denotes the Dirac delta function at~$1$.
There are similar integral formulas invol\-ving~$L_{-n}$, but then we must use an approximation to
$\delta_1(\rho_{\infty})$.

\section{Calculating moments}
\label{moments}

Throughout this section, to simplify notation, we will write $E(\cdot)=\int (\cdot)d\nu_0$.

\subsection{The basic idea}

Suppose that $n>0$.
The basic observation is that if $p(u)$ is homogeneous of degree~$n$, where $\deg(u_j)=j$, then
$
\overline{L}_{-n}(\rho_0^n p(u))
$
does not depend upon $\rho_0$.
Recall also that $\overline{L}_{-n}(\rho_{\infty})=0$.
We can now apply inf\/initesimal invariance to obtain
\begin{gather*}
E\big(\overline{L}_{-n}\big(\rho_0^n p(u)\big)\big)=\int\overline{L}_{-n}\big(\rho_0^n p(u)\big)\delta_1(\rho_{\infty})d\mu_0=0,
\end{gather*}
which gives rise to an integral formula.

To prove Theorem~\ref{2moment}, we use the identity
\begin{gather}
\label{recursion}
\overline{L}_{-1}(\rho_0u_{n+1}\bar{u}_n) =-(n+2)u_{n+1}\bar{u}_{n+1}+(n+1)u_n\bar{u}_n.
\end{gather}

\begin{Theorem}
\begin{gather*}
\int u_nu_n^*d\nu_0=\frac1{n+1}.
\end{gather*}
\end{Theorem}

\begin{proof}
Formula~\eqref{recursion}, together with inf\/initesimal invariance, implies the recursion relation
\begin{gather*}
-(n+2)E(u_{n+1} \bar{u}_{n+1})+(n+1)E(u_n\bar{u}_n)=0
\end{gather*}
with the initial condition $E(u_0 \bar{u}_0)=E(1)=1$.
\end{proof}

We will use the following notation throughout this section.

\begin{Definition}\quad
\begin{enumerate}\itemsep=0pt
\item[(a)]
\begin{gather*}
\mathbb C[u]^{(n)}:=\spanop\left\{\prod\limits_{k\ge1}u_k^{p_k}:\sum\limits_{k\ge1} kp_k=n\right\},
\end{gather*}
i.e.~the $e^{in}$ eigenspace for the action of rotations; $\mathbb C[\bar{u}]^{(n)}$ is def\/ined similarly.

\item[(b)] For each $n\geq1$ we denote
\begin{gather*}
\mathbb C[u,\bar{u}]^{(n,n)}:=\spanop\left\{\prod\limits_{k\ge1}u_k^{p_k}\bar{u}_k^{q_k}:\sum\limits_{k\ge1}
kp_k=\sum\limits_{k\ge1} kq_k=n\right\}
\end{gather*}
or in other words $\mathbb C[u,\bar{u}]^{(n,n)}\cong \mathbb C[u]^{(n)}\otimes\mathbb C[\bar{u}]^{(n)}$.
We will refer to elements in the vector space $\mathbb{C}[u,\bar{u}]^{(n,n)}$ as being of \textit{level}~$n$.
\end{enumerate}
\end{Definition}

Note that the dimension of $\mathbb C[u]^{(n)}$ is $p(n)$, the number of partitions of~$n$, hence grows very rapidly.

The rationale for the notation is the following.
The outer tensor product, $\mathcal W \times \bar{\mathcal W}$, acts on the tensor product $\mathbb C[u]\otimes \mathbb
C[\bar{u}]$.
The product of the corresponding rotation groups acts, and induces a~bigrading.
In (b) we are considering the $0$-eigenspace for the real embedded rotation group.

If $x\in \mathbb{C}[u,\bar{u}]^{(n,n')}$, then one may verify
\begin{gather*}
E(x)=e^{i(n-n')}E(x)
\end{gather*}
using the rotational invariance of Werner's measure.
Therefore, we restrict ourselves to computing integrals of elements at levels $n=1,2,\ldots$ (i.e., $n=n'$).

Suppose $1\leq m\leq n$.
In general, we can obtain integral identities by computing
\begin{gather}
\label{generalm}
\overline{L}_{-m}\circ \rho_0^m: \ \mathbb{C}[u,\bar{u}]^{(n,n-m)}\longrightarrow \bigoplus_{j=0}^m
\mathbb{C}[u,\bar{u}]^{(n-j,n-j)},
\qquad
m\geq1,
\end{gather}
and applying inf\/initesimal invariance.
As we will see in the following sections, we are particularly interested in the cases $m=1,2$.

\begin{Remark}
In~\eqref{generalm} it is necessary to restrict consideration to $ \overline{L}_{-m}$ for $m\ge1$, because we actually
need this derivative to f\/ix $\rho_{\infty}$.
Otherwise we cannot apply integration by parts to obtain integrals.
\end{Remark}

We will now give an example, where we compute the integrals for all elements of level~2.
The single equation
\[
\overline{L}_{-2}\big(\rho_0^2 u_1^2\big)=-6u_1\bar{u}_1+14u_1^2\bar{u}_1^2-8u_1^2\bar{u}_2
\]
 implies
 \[
 14
E\big(u_1^2\bar{u}_1^2\big)=6 E(u_1\bar{u}_1)+8E\big(u_1^2\bar{u}_2\big)
\] by inf\/initesimal invariance.
Therefore,
\begin{gather}
\label{level2a}
14 E\big(u_1^2\bar{u}_1^2\big)=3+8E\big(u_1^2\bar{u}_2\big)
\end{gather}
by Theorem~\ref{2moment}.
On one hand,
\[
\overline{L}_{-1}\big(\rho_0 u_1^2\bar{u}_1\big)=-3u_1^2\bar{u}_2+2u_1\bar{u}_1.
\]
On the other hand,
\[
L_{-1}\big(\rho_0 \bar{u}_1^2u_1\big)=-3\bar{u}_1^2u_2+2u_1\bar{u}_1.
\]
Therefore,
\begin{gather*}
E\big(u_1^2\bar{u}_2\big)=E\big(u_2\bar{u}_1^2\big)=\frac{1}{3}.
\end{gather*}
Equation~\eqref{level2a} can now be used to obtain the following.

\begin{Proposition}
\begin{gather*}
E(u_2\bar{u}_2)=\frac1{3},
\qquad
E\big(u_2\bar{u}_1^2\big)=E\big(\bar{u}_2u_1^2\big)=\frac13,
\qquad
E\big(\vert u_1\vert^4\big)=\frac{17}{42}.
\end{gather*}
\end{Proposition}

\subsection[Expressions for $\overline{L}_{-1}$]{Expressions for $\boldsymbol{\overline{L}_{-1}}$}

Consider~\eqref{generalm} in the case $m=1$.
The f\/irst expression we derive for this operator is purely algebraic.

\begin{Lemma}
\label{alglemma}
Suppose that $\sum kp_k=n$ and $\sum kq_k=n-1$ and let $u^p\bar{u}^q:= \prod\limits_{k}u_k^{p_k}\bar{u}_k^{q_k}$.
Then
\begin{gather*}
\overline{L}_{-1}\big(\rho_0u^p\bar{u}^q\big)=\left(\sum\limits_{j\ge1}jp_ju_j^{-1}u_{j-1}u^p\bar{u}^q\right)+\left(\bigg(2\sum\limits_{j\ge1}
q_j-2\bigg)\bar{u}_1u^p\bar{u}^q\right)\\
\hphantom{\overline{L}_{-1}\big(\rho_0u^p\bar{u}^q\big)=}{}
-\left(\sum\limits_{j\ge1}(j+2)q_j\bar{u}_j^{-1}\bar{u}_{j+1}u^p\bar{u}^q\right).
\end{gather*}
The first sum of terms are of level $n-1$, and the other terms are of level~$n$.
\end{Lemma}

\begin{proof}
We calculate
\begin{gather*}
\overline{L}_{-1}\big(\rho_0u^p\bar{u}^q\big)=-\bar{u}_1u^p\bar{u}^q+
\sum\limits_{j\ge1}\left(p_ju_j^{p_j-1}j(u_{j-1}-u_1^*u_j) \prod\limits_{k\ne
j}u_k^{p_k}\prod\limits_{k\ge1}\bar{u}_k^{q_k} \right.\\
\left.\hphantom{\overline{L}_{-1}\big(\rho_0u^p\bar{u}^q\big)=}{}
+q_j\bar{u}_j^{q_j-1}(j+2)(u_1u_j-u_{j+1})^*
\prod\limits_{k\ge1}u_k^{p_k}\prod\limits_{k\ne j}\bar{u}_k^{q_k}\right)
\\
\hphantom{\overline{L}_{-1}\big(\rho_0u^p\bar{u}^q\big)}{}
=-\bar{u}_1u^p\bar{u}^q+\sum\limits_{j\ge1}\big(jp_ju_j^{-1}u_{j-1}+((j+2)q_j-jp_j)\bar{u}_1\\
\hphantom{\overline{L}_{-1}\big(\rho_0u^p\bar{u}^q\big)=}{}
-(j+2)q_j\bar{u}_j^{-1}\bar{u}_{j+1}\big)
\prod\limits_{k}u_k^{p_k}\bar{u}_k^{q_k}.
\end{gather*}
This simplif\/ies to the expression in the statement of the lemma.
\end{proof}

The second expression is in terms of divergence-type dif\/ferential operators.
We also note that the homogeneity condition on the domains can be expressed in terms of divergence-type operators.

\begin{Proposition}\label{RandN}
Let $n\geq1$.
\begin{enumerate}\itemsep=0pt
\item[$(a)$]
\begin{gather*}
\mathbb C[u]^{(n)}=\bigg\{P\in \mathbb C[u]:\sum\limits_{j\ge1}ju_j\frac{\partial P}{\partial u_j}=n P\bigg\}.
\end{gather*}

\item[$(b)$] Suppressing $\rho_0$, the map
\begin{gather*}
\mathbb C[u]^{(n)}\otimes\mathbb C[\bar{u}]^{(n-1)}\overset{\proj\circ\overline{L}_{-1}}{\longrightarrow}
\mathbb C[u]^{(n)}\otimes\mathbb C[\bar{u}]^{(n)}
\end{gather*}
is of the form $1\otimes \overline{R}_1$, where
\begin{gather*}
\overline{R}_1=\sum\limits_{k\ge1}(2\bar{u}_1\bar{u}_{k}-(k+2)\bar{u}_{k+1})\frac{\partial}{\partial
\bar{u}_k}-2\bar{u}_1.
\end{gather*}

\item[$(b')$] The linear map $\overline{R}_1$ is injective.

\item[$(c)$] Similarly,
\begin{gather*}
\mathbb C[u]^{(n)}\otimes\mathbb C[\bar{u}]^{(n-1)}\overset{\proj\circ\overline{L}_{-1}}{\longrightarrow}
\mathbb C[u]^{(n-1)}\otimes\mathbb C[\bar{u}]^{(n-1)}
\end{gather*}
is of the form $N_1\otimes1$, where
\begin{gather*}
N_1=\sum\limits_{j\ge1}ju_{j-1}\frac{\partial}{\partial u_j}.
\end{gather*}
\end{enumerate}
\end{Proposition}

\begin{proof}
We will prove (b$'$): If $n=1$, then $\overline{R}_1:\mathbb{C}\to \mathbb{C}[\bar{u}]^{(1)}$ is injective by dimension
conside\-ra\-tions.
If $n\geq~2$, then consider the representation $\bar{\pi}$ of $\overline{\mathcal W}$ on $\mathbb C[\rho_0,u]$.
For the lowest-weight representation generated by $\rho_0^{n-1}$, we have $c=0$ and $h=-(\frac{n-1}{2})$ (see
Section~\ref{representations}).
This is a~reducible Verma module if and only if
\begin{gather*}
-(n-1)=\frac{m^2-1}{12}.
\end{gather*}
When the Verma module is irreducible, the creation operator $\overline{L}_{-1}$ is injective at each level,
i.e.~$\overline{R}_1$ is injective.
Notice that the same thing would be true for $\overline{L}_{-k}$ for any $k>0$.
\end{proof}

\begin{Remark}
For $u^P\bar{u}^Q\in\mathbb C[u]^{(n)}\otimes\mathbb C[\bar{u}]^{(n)}$ such that $\bar{u}^Q$ is in the image of
$\overline{R}_1$, there is a~recursion formula
\begin{gather*}
E\big(u^P\bar{u}^Q\big)=E\left((N_1\otimes1)\big(u^P\overline{R}_1^{-1}\big(\bar{u}^Q\big)\big)\right)=E\left(N_1\big(u^P\big)
\overline{R}_1^{-1}\big(\bar{u}^Q\big)\right),
\end{gather*}
where we are denoting a~partial inverse to $\overline{R}_1$ by $\overline{R}_1^{-1}$.
Unfortunately this does not make any sense for most~$Q$.
\end{Remark}

\begin{Definition}
For a~single complex variable~$z$, we def\/ine $\frac1{\sqrt{k!}}z^k$ to be an orthonormal basis for $\mathbb C[z]$.
For a~tensor product such as
\begin{gather*}
\mathbb C[u]=\mathbb C[u_1]\otimes \mathbb C[u_2]\otimes \cdots
\end{gather*}
we take the tensor product Hilbert space structure, meaning that $\frac{1}{\sqrt{p!}}u^p$ is an orthonormal basis, where
\begin{gather*}
p!:=p_1!p_2!\cdots.
\end{gather*}
\end{Definition}

\begin{Proposition}\quad
\begin{enumerate}\itemsep=0pt
\item[$(a)$] The adjoint of
\begin{gather*}
\mathbb C[\bar{u}]^{(n-1)}\overset{\overline{R}_1}{\longrightarrow}\mathbb C[\bar{u}]^{(n)},
\end{gather*}
where
\begin{gather*}
\overline{R}_1=\sum\limits_{k\ge1}(2\bar{u}_1\bar{u}_{k}-(k+2)\bar{u}_{k+1})\frac{\partial}{\partial
\bar{u}_k}-2\bar{u}_1
\end{gather*}
is
\begin{gather*}
\mathbb C[\bar{u}]^{(n-1)}\overset{\overline{R}_1^t}{\longleftarrow}\mathbb C[\bar{u}]^{(n)}
\end{gather*}
given~by
\begin{gather*}
\overline{R}_1^t=\sum\limits_{k\ge1}\left(2\bar{u}_k\frac{\partial}{\partial \bar{u}_1}\frac{\partial}{\partial
\bar{u}_k}-(k+2)\bar{u}_{k}\frac{\partial}{\partial \bar{u}_{k+1}}\right)-2\frac{\partial}{\partial \bar{u}_1}.
\end{gather*}

\item[$(b)$] Let~$K$ denote the kernel of $R_1^t$, i.e.~the cokernel of $R_1$ $($or the orthogonal complement of the image of
$R_1)$.
Then
\begin{gather*}
0 \longleftarrow \mathbb C[u]^{(n-1)} \overset{R_1^t}{\longleftarrow} \mathbb C[u]^{(n)} \longleftarrow K
\longleftarrow0,
\end{gather*}
i.e., $R_1^t$ is surjective.

\item[$(c)$]
\begin{gather*}
\big(\Image(1\otimes \overline{R}_1)+\Image(R_1\otimes1)\big)^{\perp}=\kernel\big(1\otimes\overline{R}_1^t\big)\cap\kernel\big(R_1^t\otimes1\big)
\\
\hphantom{\big(\Image(1\otimes \overline{R}_1)+\Image(R_1\otimes1)\big)^{\perp}}{}
=\mathbb C[u]^{(n)}\otimes \overline{K} \cap K\otimes \mathbb C[\bar{u}]^{(n)}=K\otimes \overline{K},
\end{gather*}
which has dimension $(p(n)-p(n-1))^2>0$ for $n>1$ $(p(\cdot)$ is the partition function$)$.
\end{enumerate}
\end{Proposition}

\begin{proof}
Because of the normalization for the Hermitian inner product, the adjoint for multiplication by~$z$ on $\mathbb C[z]$ is
$\frac{\partial}{\partial z}$ on the $\mathbb C[z]$, and vice versa.
This leads to the formula for $\overline{R}_1^t$.

Part (b) follows from the injectivity of $R_1$ (see (b$'$) of Proposition~\ref{RandN}).

Part (c) is elementary linear algebra: for the sum of two subspaces, the annihilators is the intersection of the
annihilators.
\end{proof}

\begin{Example}
When $n=2$,
\begin{gather*}
\kernel\big(R_1^t\big)=\mathbb C\big\{u_1^2\big\}.
\end{gather*}

When $n=3$,
\begin{gather*}
\kernel\big(R_1^t\big)=\mathbb C\big\{u_1^3+2u_1u_2\big\}.
\end{gather*}
Note $p(3)-p(2)=3-2=1$.

When $n=4$,
\begin{gather*}
\kernel\big(R_1^t\big)=\mathbb{C}\big\{4u_2^2-6u_1u_3,~3u_1^2+16u_1^2u_2+16u_1u_3\big\}.
\end{gather*}
Note $p(4)-p(3)=5-3=2$.
\end{Example}

We will now give a~slight generalization of Theorem~\ref{2moment} using the algebraic expression for~$\overline{L}_{-1}$.

\begin{Corollary}
Suppose that $\weight(p)=n$.
Then
\begin{gather*}
E\big(u^p\bar{u}_{n}\big)=\frac1{n+1}.
\end{gather*}

\end{Corollary}

\begin{proof}
The formula in Lemma~\ref{alglemma} implies
\begin{gather*}
\overline{L}_{-1}\big(\rho_0u^p\bar{u}_{n-1}\big)=\sum\limits_{j\ge1}jp_ju_j^{-1}u_{j-1}u^p\bar{u}_{n-1}-(n+1)u^p\bar{u}_n.
\end{gather*}
Thus we obtain a~recursion relation
\begin{gather*}
(n+1)E\big(u^p\bar{u}_n\big)=\sum\limits_{j\ge1}jp_jE\big(u_j^{-1}u_{j-1}u^p\bar{u}_{n-1}\big).
\end{gather*}
The terms on the right hand side of the same form with $\weight=n-1$.
Since $\sum\limits_{j\ge1}jp_j=n$, induction implies the right hand side equals~$1$.
This implies the corollary.
\end{proof}

\subsection[Expressions for ${\overline{L}_{-2}}$]{Expressions for $\boldsymbol{\overline{L}_{-2}}$}

We now consider the operator~\eqref{generalm} in the case $m=2$, which is substantially more complicated than in the
$m=1$ case.
Recall that $p_k^{(-1)}$ denotes the Laurent coef\/f\/icient of $\frac{z}{u(z)}$ and $P_2=7u_1^2-4u_2$.

\begin{Proposition}
Let $n\geq~2$.
\begin{enumerate}\itemsep=0pt
\item[$(a)$] Suppressing $\rho_0^2$, the map
\begin{gather*}
\mathbb C[u]^{(n)}\otimes\mathbb C[\bar{u}]^{(n-2)}\overset{\proj\circ
\overline{L}_{-2}}{\longrightarrow}\mathbb C[u]^{(n)}\otimes\mathbb C[\bar{u}]^{(n)}
\end{gather*}
is of the form $1\otimes \overline{R}_2$, where
\begin{gather*}
\overline{R}_2=2\bar{P}_2-\sum\limits_{j=1}^{\infty}\Big(\bar{P}_2\bar{u}_j-3(j+2)\bar{u}_1\bar{u}_{j+1}+(j+3)\bar{u}_{j+2}-\bar{p}_{j+2}^{(-1)}
\Big)\frac{\partial}{\partial \bar{u}_j}.
\end{gather*}

\item[$(a')$] The linear map $\overline{R}_2$ is injective.

\item[$(b)$] Similarly,
\begin{gather*}
\mathbb C[u]^{(n)}\otimes\mathbb C[\bar{u}]^{(n-2)}\overset{\proj\circ \overline{L}_{-2}}{\longrightarrow}
\bigoplus_{j=1}^2 \mathbb C[u,\bar{u}]^{(n-j,n-j)}
\end{gather*}
is of the form $N_2\otimes1-3N_1\otimes \bar{u}_1$, where
\begin{gather*}
N_2=\sum\limits_{j\ge~2}(j-1)u_{j-2}\frac{\partial}{\partial u_j}.
\end{gather*}

\item[$(c)$] If $u^P\otimes\bar{u}^Q\in \mathbb{C}[u,\bar{u}]^{(n,n)}$ such that $\bar{u}^Q$ lies in the image of
$\overline{R}_2$, then
\begin{gather*}
E\big(u^P\bar{u}^Q\big)=E\Big(N_2\big(u^P\big) \overline{R}_2^{-1}\big(\bar{u}^Q\big)\Big)-3E\Big(N_1\big(u^P\big)\bar{u}_1
\overline{R}_2^{-1}\big(\bar{u}^Q\big)\Big).
\end{gather*}
\end{enumerate}
\end{Proposition}

\begin{proof}
The proof of (a$'$) is the same as (b$'$) of Proposition~\ref{RandN}.
Parts (a) and (b) follow by the formulas
\begin{gather*}
\bar{\pi}(L_{-2})\big(\rho_0^2\big)=\bar{P}_2,
\qquad
\rho_0^2\bar{\pi}(L_{-2})(u_j)=\frac{j}{2}u_j\bar{P}_2-3ju_{j-1}\bar{u}_1+(j-1)u_{j-2},
\end{gather*}
and
\begin{gather*}
\rho_0^2\pi(L_{-2})\big(u(z)\big)=\frac1{u(z)}-\left(\frac1{z}-3u_1\right)u'(z)-\frac12 P_2(zu'(z)+u(z)),
\end{gather*}
which we then expand to obtain $\rho_0^2\pi(L_{-2})(u_j)$.

Applying inf\/initesimal invariance to $\overline{L}_{-2}\big(\rho_0^2 u^P\bar{u}^Q\big)$ gives part~(c).
\end{proof}

\begin{Proposition}
Fix $n\geq~2$ and let $K_m=\kernel(R_m^t: \mathbb{C}[u]^{(n)}\to \mathbb{C}[u]^{(n-m)})$ for $m=1,2$.
Then
\begin{gather*}
K_1\cap K_2=\{0\}
\end{gather*}
or
\begin{gather*}
\image(R_1)+\image(R_2)=\mathbb C[u]^{(n)}.
\end{gather*}
Therefore, in principle, we can determine all moments by using only $\overline{L}_{-1}$ and $\overline{L}_{-2}$.
\end{Proposition}

\begin{proof}
Consider the cylic~$\pi$-representation generated by $\rho_0^{n}$:
\begin{gather*}
\pi\big(\mathcal{U}(\mathcal{W})\big)\rho_0^{n}=\bigoplus_{k=0}^{\infty}\rho_0^{n-k}\mathbb{C}[u]^{(k)},
\end{gather*}
which is an irreducible Verma module.
Therefore, the~$n$-th graded component, $\mathbb{C}[u]^{(n)}$, has a~basis consisting of elements of the form
\begin{gather*}
L_{-i_j}\cdots L_{-i_1} \big(\rho_0^n\big),
\end{gather*}
where $0<i_1\leq \cdots\leq i_j$ and $i_1+\cdots+ i_j=n$.
The claim follows since $\mathcal{U}(\bigoplus_{k\geq1} \mathbb{C}L_{-k})$ is generated by $L_{-1}$ and $L_{-2}$.
\end{proof}

Consider $u^P\otimes\bar{u}^Q\in \mathbb C[u]^{(n)}\otimes \mathbb C[\bar{u}]^{(n)}$.
In principle, we can write
\begin{gather*}
\bar{u}^Q=\overline{R}_{1}(\bar{f}_1)+\overline{R}_{2}(\bar{f}_2)
\end{gather*}
for some polynomials $\bar{f}_j\in \mathbb{C}[\bar{u}]^{(n-j)}$.
We can then compute
\begin{gather*}
E\big(u^P\bar{u}^Q\big)=E\big(N_{1}\big(u^P\big)(\bar{f}_1-3\bar{u}_1\bar{f}_2)\big)+E\big(N_{2}\big(u^P\big)\bar{f}_2\big).
\end{gather*}
The question now becomes how to divide $\bar{u}^Q$ into two pieces.
In theory, this can be done using the orthogonal decomposition
\begin{gather*}
\mathbb C[\bar{u}]^{(n)}=\image(\overline{R}_{1})\oplus \left(\image(\overline{R}_{2})\ominus\image(\overline{R}_{1})\right).
\end{gather*}

\begin{Remark}
This gives a~recursion relation for moments.
The drawback is that we have to f\/ind all of the moments at a~given level (indexed by~$n$, which involves
$u_1,\dots,u_n$) to proceed to the next level.
In implementing this procedure numerically (e.g.\
for the purpose of trying to reconstruct the distribution for $u_1$), we have found it convenient to not take the
orthogonal complement, i.e.~to work with an overdetermined system of linear equations.
This has the advantage of providing consistency checks for all of our calculations.
However, because~$p(n)$ grows quite rapidly, this is slow (As of this writing, we do not have a~conjecture for the
distribution of~$u_1$).
\end{Remark}

\subsection{Uniqueness of Werner's measures}

To close this section, we will now give an alternate proof of the uniqueness of Werner's family of measures (when
$\nu_0$ is normalized to be a~probability measure).
Our statement is marginally stronger than Werner's, in that we only need to assume the measures are locally f\/inite,
i.e.~$0<\mu_0(\{a<\vert z\vert< A\})<\infty$ for some f\/inite $0<a<A<\infty$ (this is implied by the nontriviality
condition~\eqref{nontrivial}, but not vice versa).

\begin{Theorem}
If there exists a~family of locally finite measures $\{\mu_S\}$ on self-avoiding loops on Riemann surfaces which
satisfies conformal restriction, then this family is unique up to multiplication by an overall positive constant.
\end{Theorem}

\begin{proof}
We f\/irst claim that $\mu_0$ is uniquely determined (up to a~constant which we can normalize).
Conformal invariance of $\mu_0$ implies that there is a~factorization as in Proposition~\ref{technical}:
\begin{gather*}
d\mu_0(\gamma)=d\nu_0(u) \times \frac{d\rho_{\infty}}{\rho_{\infty}},
\end{gather*}
where here we view $d\nu_0$ as a~measure on~$u$.
Local f\/initeness implies that $\nu_0$ is f\/inite (see the proof of (b) of Proposition~\ref{introlemma}), and hence we can
normalize it to be a~probability measure.
The measure $\nu_0$ is completely determined by the joint distributions of $u_1,\dots,u_N$, $N\ge1$.
The~$u_j$ are bounded, hence these distributions are determined by their joint moments.
Finally we have shown, using only inf\/initesimal conformal invariance (which depends on the conformal restriction
property, and not any specif\/ic features of Werner's construction), that the moments for these distributions are (in
principle) computable.
This determines~$\mu_0$.

The proof that $\mu_0$ determines $\mu_S$ for all~$S$ basically follows from the argument given in Section~6.1
of~\cite{W}.
However there is a~slight f\/law in that argument.
It is not quite the case that ``The family of events of the type
\begin{gather*}
A_D=\{\gamma:
\
\gamma\subset D~\text{and goes around the hole in}~D\},
\end{gather*}
when~$D$ varies in the family of annular regions in~$S$ is stable under f\/inite intersections''.
For example in the plane the intersection of the two annuli $A_1:=\Delta\setminus\{\vert z\vert<1/8\}$ and
$A_2:=\Delta\setminus\{\vert z-1/2\vert<1/8\}$ is a~pair of pants; there does not exist an annulus inside of this pair
of pants which contains all the loops which go around both holes.
So the argument must be modif\/ied (this kind of argument is also used in the earlier proof of Lemma~4 of~\cite{W}, and in
that context it is valid, because $0$ is always assumed to be in the hole of the allowed annuli).

Given knowledge of $\mu_0$, for any proper open subset~$S$ of the plane, and for any nontrivial free homotopy class
$C\subset \Loop(S)$, $\mu_S(C)$ is uniquely determined; this follows from conformal restriction, because we can
assume $S\subset \mathbb C\setminus\{0\}$ and all the loops in~$C$ go around zero.

Suppose~$S$ is a~general Riemann surface.
Consider the family of events $C=C_D$, where~$C$ is a~nontrivial free homotopy class of loops in an open subset
$D\subset S$ such that~$D$ is conformally equivalent to a~proper open subset of $\mathbb C$.
We claim this family is stable under f\/inite intersections.
Clearly $D_1\cap D_2$ is conformally equivalent to a~proper open subset of $\mathbb C$.
The main point is to show that $C_1\cap C_2$ determines a~unique free homotopy class in $D_1\cap D_2$.
This topological fact is probably well-known, but we will give a~proof.

Suppose that we are given a~f\/ixed conformal equivalence of~$D$ with a~proper open subset of~$\mathbb C$, $i:D\to \mathbb
C$, and a~free homotopy class $C\subset \Loop(D)$.
If $\gamma\in \Loop(D)$, then (by the Jordan curve theorem applied to~$i(\gamma)$) the complement of $i(D)$ is
divided into an inside, $\In_D(\gamma)$, and an outside, $\Out_D(\gamma)$ (which contains~$\infty$, i.e.~large~$z$).

\begin{Lemma}
If $\gamma_1\in C_D$, then $C_D$ is determined by $\In_D(\gamma_1)$, i.e.~if $\gamma_2\in Loop(D)$, then
$\gamma_2\in C_D$ if and only if $\In_D(\gamma_2)=\In_D(\gamma_1)$.
\end{Lemma}

\begin{proof}
This is a~topological claim, so in a~standard way we can suppose loops are smooth, and intersections are transverse.
Suppose $H(s,t)$ is a~homotopy (with $\Image(H(\cdot,0))=\gamma_1$, and $\Image(H(\cdot,1))=\gamma_2$).
For $0<t<1$, $H(\cdot,t)$ is not necessarily simple, but we can nonetheless talk about $\In_D(H(\cdot,t))$, by
using the inner boundary.
This set, $\In_D(H(\cdot,t))$, is independent of~$t$, by continuity, and this implies $\In_D(\gamma_2)=\In_D(\gamma_1)$.

Now consider the converse.
Let $U_{\pm}^{(j))}$ denote the bounded and unbounded components for $\mathbb C\setminus \gamma_j$, respectively.
Then $U_+^{(1)}\cap U_+^{(2)}$ and $\{\infty\}\cup (U_-^{(1)}\cap U_+^{(2)})$ are open contractible sets (for
example $U_+^{(1)}\cap U_+^{(2)}$ is the bounded component for the inner boundary of $\gamma_1 \cup\gamma_2$).
The complement is a~closed region with a~boundary composed of the inner and outer boundaries for $\gamma_1
\cup\gamma_2$, and it is homotopic to a~annulus with boundary (for curves which intersect transversely, it is an annulus
which is pinched at the points of intersection of the $\gamma_i$).
This annular region is entirely contained in~$D$, and hence $\gamma_1$ and $\gamma_2$ are homotopic in~$D$.
\end{proof}

We now use this to show that $C_1\cap C_2$ determines a~unique free homotopy class in $D_1\cap D_2$.
Fix conformal embeddings $i_j:D_j\to\mathbb C$, and use the restriction of $i_1$ to embed $D_1\cap D_2$.
Suppose that $\gamma_1,\gamma_2\in C_1\cap C_2$.
Then the Lemma implies that $\In_{D_k}(\gamma_j)$ does not depend on~$j$.
But then $\In_{D_1\cap D_2}(\gamma_j)$ is also independent of~$j$, and hence by the Lemma, $\gamma_1$ and
$\gamma_2$ are homotopic in $D_1\cap D_2$.

This now implies that $C_1\cap C_2$ determines a~unique free homotopy class in $D_1\cap D_2$.
This class is clearly nontrivial, because its image in $D_k$ is nontrivial, $k=1,2$.
This now implies that the set of events $C_D$ is stable under f\/inite intersections.
Now the argument in Section~6.1 of~\cite{W} implies $\mu_S$ is uniquely determined.
\end{proof}

\section{The diagonal distribution}
\label{diagonaldistribution}

To determine the joint distribution for $(\rho_0,\rho_{\infty})$, Proposition~\ref{introlemma} implies that it suf\/f\/ices
to determine the distribution for $H=-\log (a)\ge0$, which by part~(e) of Proposition~\ref{introlemma} is a~kind
of height function for
\begin{gather*}
\big\{\sigma\in \Homeo\big(S^1\big):\exists \  \text{unique welding}~\sigma=lau\big\}.
\end{gather*}

\begin{Conjecture}
For some $\beta_0<\frac{5\pi^2}{4}$, the $\nu_0$ distribution for~$a$ is given~by
\begin{gather*}
\nu_0(\{\sigma:\exp(-x)\le a(\sigma)\le1\})= \exp(-\beta_0/x),
\qquad
x>0.
\end{gather*}
Equivalently the Laplace transform
\begin{gather}
\label{laplacetransform}
\int a^{\lambda} d\nu_0(\sigma)
=\int_{x=0}^{\infty} a^{\lambda} d\exp(-\beta_0/x)=2\sqrt{\lambda\beta_0}K_1\big(2\sqrt{\lambda \beta_0}\big)
\end{gather}
for $\lambda>0$, where $K_1$ is a~modified Bessel function.
\end{Conjecture}

We will f\/irst explain how this conjecture is related to a~remarkable calculation of Werner in Section~7 of~\cite{W}.
We will then present some calculations which are possibly relevant to a~proof, and incidentally give an estimate for
Werner's constant.
Finally we will brief\/ly indicate how the conjecture naturally generalizes to the deformation of Werner's measure
considered in~\cite{KS}.

\subsection{A formula of Werner}

As in Section~7 of~\cite{W}, consider the function
\begin{gather*}
F(\rho):=\mu\big(\Loop^1(A)\big),
\end{gather*}
where~$A$ is a~f\/inite type annulus with modulus $\rho=\rho(A)$, i.e.~$\rho>0$ is the unique number such that~$A$ is
conformally equivalent to
\begin{gather*}
\{1<\vert z\vert<e^{\rho}\}.
\end{gather*}

Cardy (see~\cite{Cardy}) has conjectured an exact formula
\begin{gather}
\label{cardyformula}
F(\rho)=6\pi\frac{\sum\limits_{k\in\mathbb Z}(-1)^{k-1}kq^{3k^2/2-k+1/8}}{\prod\limits_{k=1}^{\infty}(1-q^k)},
\qquad
q=\exp\big({-}2\pi^2/\rho\big).
\end{gather}
As we will explain below in more detail
\begin{gather}
\label{setinclusions}
\Loop^1\big(\big\{1<\vert z\vert<e^{\rho}\big\}\big)\subset\big\{1\le \rho_0\le \rho_{\infty}\le e^{\rho}\big\}\subset
\Loop^1\left(\left\{\frac14<\vert z\vert<4e^{\rho}\right\}\right)
\end{gather}
and as a~consequence
\begin{gather}
\label{intertwine}
F(\rho)\le \int_0^{\rho}\nu_0\big(e^{-x}\le a\le1\big)dx\le F(\log (16)+\rho).
\end{gather}
This incidentally explains the constant $6\pi$ in~\eqref{cardyformula}, which ensures that the derivative of~$F$ is
asymptotically one, or equivalently that $\nu_0$ is a~probability measure.

Werner shows that $F(\rho)$ is asymptotic to ${\rm const} \cdot \exp \big({-}\frac{\beta}{\rho}\big)$ as $\rho\to0$, where
$\beta=\frac{5\pi^2}{4}$; see Proposition~18 of~\cite{W}.
This leads to the upper bound on $\beta_0$ in our statement of the diagonal distribution conjecture (if Cardy's
conjecture is correct, then~\eqref{intertwine} implies sharper upper and lower bounds for $\beta_0$).

\begin{Lemma}\label{lemma7.2}
Fix $x>0$.
If $\gamma\in \Loop^1(\{1<\vert z\vert<e^x\})$, then
\begin{gather*}
1\le \rho_0(\gamma)\le \rho_{\infty}(\gamma)\le e^x.
\end{gather*}
\end{Lemma}

\begin{proof}
The Cauchy integral formula implies, for suf\/f\/iciently smooth~$\gamma$,
\begin{gather*}
\frac1{\rho_0}=\big(\phi_+^{-1}\big)'(0)=\frac{1!}{2\pi i}\int_{\gamma}\frac{\phi_+^{-1}(t)}{t^2}dt.
\end{gather*}
Since~$\gamma$ is outside the unit disk and $\phi_+^{-1}:U_+\to \Delta$.
This implies
\begin{gather*}
\frac1{\rho_0}\le \frac{\text{length}(\gamma)}{2\pi}\le1 .
\end{gather*}
This implies the f\/irst inequality.
The last inequality also follows from this.

We noted previously that the equality
\begin{gather*}
a^2=\frac{\rho_0^2}{\rho_{\infty}^2}=\frac{1-\sum\limits_{m=1}^{\infty}(m-1)\vert
b_m\vert^2}{1+\sum\limits_{n=1}^{\infty}(n+1)\vert u_n\vert^2}
\end{gather*}
implies $a\le1$, i.e.~$\rho_0\le \rho_{\infty}$.
\end{proof}

\begin{Lemma}\label{lemma7.3}\quad
\begin{enumerate}\itemsep=0pt
\item[$(a)$]
\begin{gather*}
\mu_0\big\{\gamma:1\le \rho_0(\gamma)\le \rho_{\infty}(\gamma)\le e^x\big\}=\int_{y=0}^x\nu_0\big\{e^{-y}\le a\le1\big\}dy \le
x\nu_0\big\{e^{-x}\le a\le1\big\}.
\end{gather*}

\item[$(b)$]
\begin{gather*}
\int_0^x e^{-\beta_0/y}dy=xe^{-\beta_0/x}-\Gamma\left(0,\frac{\beta_0}{x}\right)=xe^{-\beta_0/x}-\text{Ei}\left(1,\frac{\beta_0}{x}\right)
\\
\hphantom{\int_0^x e^{-\beta_0/y}dy}{}
=x-\beta_0 \log (x)+\beta_0(\log (\beta_0)+\gamma-1)-\frac12\beta_0^2 x^{-1}+\cdots.
\end{gather*}

\item[$(c)$] There is an asymptotic expansion
\begin{gather*}
\int_0^x e^{-\beta_0/y}dy=\frac{x^2}{\beta_0}e^{-\beta_0/x}\sum\limits_{n=0}^{\infty}(-1)^n
n!\left(\frac{x}{\beta_0}\right)^n
\\
\hphantom{\int_0^x e^{-\beta_0/y}dy}{}
=\frac{x^2}{\beta_0}e^{-\beta_0/x}\left(1-\frac{x}{\beta_0}+2\frac{x^2}{\beta_0^2}-\cdots\right)
\qquad\text{as}\quad x\to0.
\end{gather*}
\end{enumerate}
\end{Lemma}

\begin{proof}
(a) Using the factorization $d\mu_0=\frac{d\rho_{\infty}}{\rho_{\infty}}\times d\nu_0$,
\begin{gather*}
\mu_0\big\{\gamma:1\le \rho_0(\gamma)\le \rho_{\infty}(\gamma)\le e^x\big\}=\mu_0\left\{\frac1{\rho_{\infty}}\le a\le1,
\
1\le \rho_{\infty}\le e^x\right\}
\\
\hphantom{\mu_0\big\{\gamma:1\le \rho_0(\gamma)\le \rho_{\infty}(\gamma)\le e^x\big\}}{}
=\int_{\rho_{\infty}=1}^{e^x}\nu_0\left\{\frac1{\rho_{\infty}}\le a\le1\right\}\frac{d\rho_{\infty}}{\rho_{\infty}}.
\end{gather*}
By making the change of variables $\rho_{\infty}=e^y$, we obtain the expression in part (a).

(b) and (c) are standard facts.
For example there is a~Laurent expansion
\begin{gather*}
e^{-\beta_0/x}=1-\frac{\beta_0}{x}+\frac12\left(\frac{\beta_0}{x}\right)^2-\cdots,
\qquad
0<\vert x\vert<\infty.
\end{gather*}
Therefore there is an expansion
\begin{gather*}
\int_0^x e^{-\beta_0/y}dy=x-\beta_0 \log (x)+c_0-\frac12\beta_0^2x^{-1} +\frac{\beta_0^3}{3!2}x^{-2}-\cdots\\
\hphantom{\int_0^x e^{-\beta_0/y}dy=}{}
-(-1)^n\frac1{n!(n-1)}\beta_0^{n}\frac1{x^{n-1}}-\cdots,
\end{gather*}
where the divergence of the logarithm and the Laurent expansion at $x=0$ perfectly cancel, allowing us to f\/igure
out~$c_0$.
\end{proof}

\begin{Corollary}
\begin{gather*}
\mu_0 \big(\Loop^1\big(\big\{1<\vert z\vert<e^{\rho}\big\}\big)\big)\le \rho\nu_0\big\{e^{-\rho}\le a\le1\big\}.
\end{gather*}
\end{Corollary}

\begin{proof}
This follows from Lemma~\ref{lemma7.2} and (a) of Lemma~\ref{lemma7.3}.
\end{proof}

Here is another approach, although not quite as sharp:
\begin{gather*}
\Loop^1\big(\big\{1<\vert z\vert<e^{\rho}\big\}\big)\subset \Loop^1\big(\big\{\vert z\vert<e^{\rho}\big\}\big)\setminus \Loop(\Delta).
\end{gather*}
Werner's formula for the measure of the latter set is $c_{\rm W}\rho$, where $c_{\rm W}$ is Werner's constant (see below).

\subsection{Werner's constant}
\label{Wernerconstant}

Recall that we have normalized Werner's family of measures by assuming that $\nu_0$ is a~probability measure.
We let $c_{\rm W}$ denote the constant such that if~$\gamma$ is a~loop which surrounds~$\Delta$,
\begin{gather*}
\mu_0\big(\Loop^1(U_+\setminus\{0\})\setminus \Loop(\Delta)\big)=c_{\rm W} \log (\rho_0(\gamma)).
\end{gather*}

\begin{Proposition}
$ c_{\rm W} \ge1$.
\end{Proposition}

\begin{proof}
On the one hand
\begin{gather*}
\Loop^1\big(\big\{1<\vert z\vert<e^{\rho}\big\}\big)\subset \Loop^1\big(\big\{\vert z\vert<e^{\rho}\big\}\big)\setminus \Loop(\Delta).
\end{gather*}
Therefore by Werner's formula for the measure of the latter set,
\begin{gather*}
F(\rho)\le c_{\rm W} \rho.
\end{gather*}
On the other hand
\begin{gather*}
\Loop^1\big(\big\{1<\vert z\vert<e^{\rho}\big\}\big)\subset\big\{1\le \rho_0\le \rho_{\infty}\le e^{\rho}\big\}\subset\Loop^1\left(\left\{\frac14<\vert z\vert<4e^{\rho}\right\}\right),
\end{gather*}
where the last inclusion uses Koebe's quarter theorem.
Therefore
\begin{gather*}
F(\rho)\le \int_0^{\rho}\nu_0\big(e^{-x}\le a\le1\big)dx\le F(\log (16)+\rho).
\end{gather*}

Because
\begin{gather*}
\nu_0\big(e^{-x}\le a\le1\big)\uparrow1 \qquad\text{as}\quad x\uparrow1
\end{gather*}
it follows that $F(\rho)$ behaves like a~linear function with slope one for $\rho\gg 1$.
This behavior is compatible with the estimate above using Werner's formula if and only if $c_{\rm W} \ge1$.
This implies the proposition.
\end{proof}

\subsection{Some ideas}

The conjectural Laplace transform~\eqref{laplacetransform} satisf\/ies the ODE
\begin{gather*}
\lambda f''(\lambda)-\beta_0 f(\lambda)=0,
\end{gather*}
Thus we need to show that
\begin{gather*}
\int\big(\lambda \log (a)^2-\beta_0\big)a^{\lambda}d\nu_0(\sigma)=0,
\qquad
\lambda>0
\end{gather*}
for some constant $\beta_0$.
Roughly speaking, we are trying to calculate the second moment for the distribution of $H=-\log (a)$.
To calculate the second moment for a~standard normal complex variable, one can apply $\partial\bar{\partial}$ to
$\exp (-\vert z\vert^2/2)$ and use inf\/initesimal invariance of the background Lebesgue measure; our strategy is to
do the same with the stress tensor $T(t)$ in place of $\partial$, $a^{\lambda}$ in place of the Gaussian, and Werner's
measure in place of Lebesgue measure.

We will now list a~number of formulas which are hopefully useful.

\begin{Lemma}\quad
\begin{enumerate}\itemsep=0pt
\item[$(a)$] For $n>0$
\begin{gather*}
L_nL_{-n}a^{\lambda}=L_{-n}L_na^{\lambda}=\frac{\lambda^2}{4}P_n(l_1,\dots,l_n)P_n(u_1,\dots,u_n)a^{\lambda-n}- n\lambda
a^{\lambda}.
\end{gather*}

\item[$(b)$] For $m>n\ge0$
\begin{gather*}
L_mL_{-n}a^{\lambda}=\frac{\lambda^2}{4}P_m(l_1,\dots,l_m)P_n(u_1,\dots,u_n)a^{\lambda}\frac{\rho_{\infty}^m}{\rho_0^n}.
\end{gather*}

\item[$(c)$]
\begin{gather*}
\overset{\longrightarrow}{L}_n\overset{\longrightarrow}{L}_{-n}a^{\lambda}
=\overset{\longrightarrow}{L}_{-n}\overset{\longrightarrow}{L}_n a^{\lambda}=
\lambda^2\Reop(P_n(l_1,\dots,l_n))\Reop(P_n(u_1,\dots,u_n))a^{\lambda-n}-~2n\lambda a^{\lambda}
\\
\hphantom{\overset{\longrightarrow}{L}_n\overset{\longrightarrow}{L}_{-n}a^{\lambda}}{}
=\lambda^2\Reop\big(P_n\big(u_1\big(\sigma^{-1}\big),\dots,u_n\big(\sigma^{-1}\big)\big)\big)\Reop(P_n(u_1,\dots,u_n))a^{\lambda-n}-~2n\lambda
a^{\lambda}.
\end{gather*}

\item[$(d)$] For $m>n\ge0$
\begin{gather*}
\overset{\longrightarrow}{L}_m\overset{\longrightarrow}{L}_{-n}a^{\lambda}=
\lambda^2\Reop(P_m(l_1,\dots,l_m))\Reop(P_n(u_1,\dots,u_n))a^{\lambda}\frac{\rho_{\infty}^m}{\rho_0^n}.
\end{gather*}
\end{enumerate}
\end{Lemma}

\begin{proof}
(a) The fact that $L_n$ and $L_{-n}$ commute when acting on $a^{\lambda}$ follows from the fact that $L_0a=0$.

Using $L_{-n}(\rho_{\infty})=0$ and (d) of Proposition~\ref{variation1},
\begin{gather*}
L_{-n}a^{\lambda}=\lambda a^{\lambda-1}L_{-n}(\rho_0)\frac1{\rho_{\infty}}=\frac{\lambda}{2\rho_0^n}P_n(u))a^{\lambda},
\end{gather*}
where we have abbreviated $P_n(u_1,\dots,u_n)=P_n(u)$.
Therefore
\begin{gather*}
L_nL_{-n}a^{\lambda}= \frac{\lambda}{2\rho_0^n} \left(\lambda a^{\lambda-1}\rho_0
L_n\left(\frac1{\rho_{\infty}}\right)P_n(u_1,\dots,u_n)+ a^{\lambda}L_n(P_n(u_1,\dots,u_n))\right).
\end{gather*}
Recall that $P_n(u_1,\dots,u_n)=-2nu_n+\text{function}(u_1,\dots,u_{n-1})$ and $L_n(u_n)=\rho_0^n$.
This implies
\begin{gather*}
L_nL_{-n}a^{\lambda}=\frac{\lambda}{2\rho_0^n} \left(\lambda a^{\lambda-1}\rho_0
\frac12(\frac1{\rho_{\infty}})^{-n+1}P_n(l)P_n(u)+ a^{\lambda}(-2n\rho_0^n)\right)
\end{gather*}
This simplif\/ies to (a).

(b) This follows in a~similar way, using the fact that $L_m$ kills $P_n(u_1,\dots,u_n)$.

(c) and (d) are proven in a~similar way, and will not be used.
\end{proof}

Recall that
\begin{gather*}
\big(\partial \log \big(\phi_+^{-1}\big)\big)^2=\left(\sum\limits_{n=0}^{\infty}P_n(\phi_+)t^n\right)\left(\frac{dt}{t}\right)^2
\end{gather*}
(this is a~holomorphic quadratic dif\/ferential which is well-def\/ined in $U_+$) and
\begin{gather*}
\big(\partial \log \big(\phi_-^{-1}\big)\big)^2=\left(\sum\limits_{n=0}^{\infty}P_n(\phi_-)t^{-n}\right)\left(\frac{dt}{t}\right)^2
\end{gather*}
(this is a~holomorphic quadratic dif\/ferential which is well-def\/ined in $U_-$; note that
\begin{gather*}
\left(\frac{dt}{t}\right)^2=\left(\frac{dt^{-1}}{t^{-1}}\right)^2.
\end{gather*}
The fact that these two quadratic dif\/ferentials do not have a~common domain, or at the very best, are possibly def\/ined
on the rough loop~$\gamma$, is a~crucial point.

\begin{Proposition}

\begin{gather*}
E\big(\big(\lambda P_n(\phi_+)P_n(\phi_-)-4n\big)a^{\lambda}\big)=E\big(\big(\lambda P_n(u)P_n(l)a^{-n}-4n\big)a^{\lambda}\big)=0.
\end{gather*}
\end{Proposition}

\begin{proof}
This follows from the Lemma and inf\/initesimal conformal invariance.
\end{proof}

The basic question now is whether there is a~constant $\beta_0$ such that $\lambda \log (a)^2-\beta_0$ is a~limit,
in an appropriate measure theoretic sense relative to $\nu_0$, of linear combinations of the functions $\lambda
P_n(\phi_+)P_n(\phi_-)-~2n$, as~$n$ varies.

\begin{Question}
Do there exist constants $c_n$ such that
\begin{gather*}
\sum\limits_{n=1}^{N}c_nP_n(u)P_n(l)a^{-n} \to \log (a)^2 \qquad\text{as}\quad N\to\infty
\end{gather*}
in some measure-theoretic sense relative to $\nu_0$?
\end{Question}

This is def\/initely false for all~$\sigma$.
To see this, suppose that
\begin{gather*}
\sigma=\phi_N(w_N,z)=z\frac{\big(1+\overline{w_N} z^{-N}\big)^{1/N}}{\big(1+w_N z^N\big)^{1/N}}.
\end{gather*}
In this case
\begin{gather*}
u(z)=z\big(1+w_Nz^N\big)^{-1/N},
\qquad
\frac {U'(t)}{U(t)}=\frac1{t}\frac1{1-w_N t^N}
\end{gather*}
and
\begin{gather*}
(\partial \log  U(t))^2=\big(1+2w_N t^N+3\big(w_Nt^N\big)^2+4\big(w_Nt^N\big)^3+\cdots\big)\left(\frac{dt}{t}\right)^2.
\end{gather*}
Thus for this particular~$u$
\begin{gather*}
P_n(u)=(m+1)w_N^m=-P_n(l)^*,
\qquad
n=mN
\end{gather*}
and zero otherwise.
Also
\begin{gather*}
l(t)=t\big(1+\bar{w}_Nt^{-N}\big),
\qquad
\partial \log (l(t))=\frac1{t\big(1+\bar{w}_Nt^{-N}\big)}dt=\frac{t^{N-1}}{t^N+\bar{w}_N}dt,
\\
a=\big(1-\vert w_N\vert^2\big)^{1/N},
\end{gather*}
so that
\begin{gather*}
\log (a)^2=\frac1{N^2}\log \big(1-\vert w_N\vert^2\big)^2.
\end{gather*}

If we actually have an identity, then for each $N=1,2,\dots $
\begin{gather*}
\frac1{N^2}\log \big(1-\vert w_N\vert^2\big)^2=\sum\limits_{m=0}^{\infty}c_{mN} \frac{(m+1)^2\vert w_N \vert^{2m}}{\big(1-\vert
w_N\vert^2\big)^m}.
\end{gather*}
If we set $x=\vert w_N\vert^2$, then this is equivalent to
\begin{gather*}
\log (1-x)^2=\sum\limits_{m=0}^{\infty}N^2c_{mN} (m+1)^2\left(\frac{x}{1-x}\right)^m.
\end{gather*}
This is clearly impossible: we cannot consistently solve for the constants.
Furthermore the radius of convergence for the l.h.s.\ is~$1$, and the radius of convergence for the r.h.s.\ is $\frac12$.

A more promising approach seems to be to use the stress-energy tensor.
Here is one heuristic calculation:
\begin{gather*}
E\big(T(t)T(s)a^{\lambda}\big)=\sum\limits_{n,m} E\big(L_{-n}L_{-m}a^{\lambda}\big)t^n
s^m\left(\frac{dt}{t}\right)^2\left(\frac{ds}{s}\right)^2
\\
\qquad{}
=\sum\limits_{n} E\big(L_{-n}L_{n}a^{\lambda}\big)t^ns^{-n}\left(\frac{dt}{t}\right)^2\left(\frac{ds}{s}\right)^2
\\
\qquad\quad{}\times
\lambda \sum\limits_{n} E\big((\lambda
P_n(\phi_+)P_n(\phi_-)-n)a^{\lambda}\big)t^ns^{-n}\left(\frac{dt}{t}\right)^2\left(\frac{ds}{s}\right)^2
\\
\qquad{}
=\lambda \left(\lambda E\big(\big(\partial \log \big(\phi_+^{-1}(t)\big)\big)^2\big(\partial
\log \big(\phi_+^{-1}(s)\big)\big)^2a^{\lambda}\big) -\sum\limits_{n}
n\left(\frac{t}{s}\right)^n\left(\frac{dt}{t}\right)^2\left(\frac{ds}{s}\right)^2E\big(a^{\lambda}\big)\right)\!
\\
\qquad{}
=\lambda \left(\lambda E\big(\big(\partial \log \big(\phi_+^{-1}(t)\big)\big)^2\big(\partial
\log \big(\phi_+^{-1}(s)\big)\big)^2a^{\lambda}\big) -\delta'\left(\frac{t}{s}\right)
\left(\frac{dt}{t}\right)^2\left(\frac{ds}{s}\right)^2E\big(a^{\lambda}\big)\right).
\end{gather*}
We now need to apply some kind of pairing for quadratic dif\/ferentials, and we are stymied at this point.

\subsection{KS conjecture and diagonal distribution}
\label{KSconjecture}

In~\cite{KS} Kontsevich and Suhov show that for each Riemann surface, there exists a~continuous positive determinant
line bundle $\Det\to \Loop(S)$, and these line bundles have a~natural restriction property.
They conjecture that for each ``central charge''~$c$ (in some range), there exists a~family of measures~$\mu_S$ having
values in the positive line bundle $\Det^c$ and satisfying a~conformal restriction property.
In the case $c=0$, this family is the family of measures constructed by Werner.

There is a~canonical trivialization of the determinant line bundle in genus zero, so that the conjectured KS measure can
be viewed as a~scalar measure which is invariant with respect to global conformal transformations; see Section~2.5
of~\cite{KS}.
We denote this measure restricted to $\Loop^1(\mathbb C \setminus\{0\})$ by $\mu_c$; properly normalized, this is
the Werner measure when $c=0$.

\begin{Lemma}
Assume that $\mu_c$ exists.
Then
\begin{enumerate}\itemsep=0pt
\item[$(a)$] The distributions for $\rho_0$ and $\rho_{\infty}$ are scale invariant.

\item[$(b)$]
\begin{gather*}
d(W_*\mu)(\sigma,\rho_{\infty})=d\nu_c(\sigma) \times \frac{d\rho_{\infty}}{\rho_{\infty}}.
\end{gather*}

\item[$(c)$] The measure $d\nu_c(\sigma)$ is inversion invariant and invariant with respect to conjugation by $C:z\mapsto z^*$.

\item[$(d)$] The measure $d\nu_c(\sigma)$ is supported on~$\sigma$ having triangular factorization $\sigma=lau$, i.e.~$m=1$.

\item[$(e)$] If in addition $\nu_c$ is finite, and hence can be normalized to be a~probability measure, then there is an
inequality generalizing~\eqref{intertwine},
\begin{gather*}
F_c(\rho)\le \int_0^{\rho}\nu_c\big(e^{-x}\le a\le1\big)dx\le F_c(\log (16)+\rho),
\end{gather*}
where $F_c(\rho):=\mu_c(\Loop^1(A))$, where~$A$ is a~finite type annulus with modulus $\rho=\rho(A)$.
\end{enumerate}
\end{Lemma}

This is a~rigorous result (contingent on the existence of $\mu_c$), because (a)--(d) use only global conformal invariance
of $\mu_c$, and (e) only depends on~\eqref{setinclusions}.

There is a~natural conjecture for the diagonal distribution (there may be a~conjecture for~$F_c(\rho)$ which is implicit
in~\cite{Cardy}, but we will not pursue this).

\begin{Conjecture}
The $\nu_c$ distribution for $H=-\log (a)$ is the inverse gamma distribution with parameters $\alpha=1-c$ and some
$\beta_c>0$ $($possibly proportional to $h^+(c)$, the larger value of two values of the conformal anomaly~$h$
corresponding to $c<1)$.
In other words we are conjecturing that
\begin{gather*}
\nu_c(\{\sigma:\exp (-x)\le a(\sigma)\le1\})=\frac{\Gamma(\alpha,\beta_c/x)}{\Gamma(\alpha)},
\qquad
x>0
\end{gather*}
and the Laplace transform
\begin{gather*}
\int a^{\lambda}d\nu_c(\sigma)=\frac{2(\beta_c \lambda)^{\frac{\alpha}{2}}}{\Gamma(\alpha)}K_{\alpha}\big(\sqrt{4\beta_c\lambda}\big),
\end{gather*}
where $K_{\alpha}$ is a~modified Bessel function.
This function of~$\lambda$ satisfies the differential equation
\begin{gather*}
\lambda f''(\lambda)+c f'(\lambda)-\beta_c f(\lambda)=0.
\end{gather*}
\end{Conjecture}

This dif\/ferential equation obviously makes sense for values of the parameters which are not necessarily positive.
But for example if $c=1$, i.e.~$\alpha=0$, then the particular solution we are considering, $K_0$, is not f\/inite at
$\lambda=0$, so that the probabilistic interpretation is lost (this is obvious by noting that the pdf is not integrable
at~$\infty$ when $\alpha=0$).
In terms of our conjecture this means that when $c=1$, the~$\sigma$ distribution for the conjectured Kontsevich--Suhov
measure is not f\/inite, according to us.

To motivate this, in a~heuristic way, we imagine that $\mu_c$ is absolutely continuous with respect to Werner's measure
$\mu_0$: $\mu_c=\delta_c d\mu_0$.
We then apply inf\/initesimal invariance in the following way.
Suppose that $n>0$.
Then
\begin{gather*}
L_n\big(L_{-n}\big(a^{\lambda}\big)\delta^c \delta(\rho_0=1)d\mu_0 \big)
=\big(\big(L_nL_{-n}\big(a^{\lambda}\big)\big)\delta^c +L_{-n}\big(a^{\lambda}\big)L_n\big(\delta^c\big)\big)\delta(\rho_0=1)d\mu_0
\\
\qquad{}
=\big(\big(\lambda^2P_n(u)P_n(l)a^{-n}-2n\lambda\big)+\lambda P_n(l)\rho_{\infty}^n
cQ_n(u,l)\big)\delta(\rho_0=1)a^{\lambda}\delta^cd\mu_c,
\end{gather*}
where we have tentatively written
\begin{gather*}
L_n(\delta_c)=cQ_n(u,l)\delta_c.
\end{gather*}
This should rigorously be expressed in terms of divergences, as proposed in Section~2.5.2 of~\cite{KS}.
From this, by dividing by~$\lambda$, we can deduce that
\begin{gather*}
\int\big(\lambda P_n(u)P_n(l)a^{-n}+cP_n(l)Q_n(u,l)\rho_{\infty}^n-2n\big)a^{\lambda}d\nu_c=0.
\end{gather*}
Now we would have to take linear combinations and limits, to obtain $\log (a)^2$ from the f\/irst term,
$\log (a)$ from the second term (involving~$c$), and a~constant $\beta_c$ from the third term.

\appendix

\section{The Vietoris topology}
\label{appendixA}

Suppose that~$S$ is a~topological space.
The Vietoris topology on $\Comp(S)$ has a~base consisting of sets of the form
\begin{gather*}
\big\{K\in \Comp(S):K\subset U, K\cap U_i\ne \phi, i=1,\dots,n\big\},
\end{gather*}
where $U,U_1,\dots,U_n$ are open subsets of~$S$.
Given $K_0\in \Comp(S)$, suppose we tightly cover $K_0$ with open sets $U_i$, $1\le i\le n$, and let
$U=\cup_iU_i$.
Then ``$K$ is close to $K_0$'' means that (i) $K\subset U$, so every point in~$K$ is close to a~point in $K_0$, and (ii)
for each point $x_0\in K_0$, $x_0\in U_i$, for some~$i$, hence $K\cap U_i \ne \varnothing$ implies $x_0$ is close to
some point in~$K$.
If~$S$ is metrizable, with metric~$d$, then the Vietoris topology is compatible with the associated Hausdorf\/f metric
topology on $\Comp(S)$, where the Hausdorf\/f metric is given by
\begin{gather*}
\delta(K_1,K_2)=\max\Big\{\sup_{p_1\in K_1}(d(p_1,K_2)),\sup_{p_2\in K_2}(d(K_1,p_2))\Big\}.
\end{gather*}
For most topological properties~$\tau$, ``$S$ is~$\tau$'' if and only if ``$\Comp(S)$ is~$\tau$'' (see Section~4
of~\cite{Michael}).
In particular if~$S$ is second countable and locally compact, then $\Comp(S)$ is second countable and locally compact.

Suppose that~$S$ is a~Riemann surface with a~f\/ixed compatible complete metric.
The associated Hausdorf\/f metric on $\Loop(S)$ is obviously not complete, since for example a~small circle can
pinch down to a~point.
Does there exist a~complete separable metric on $\Loop(S)$ compatible with the Vietoris topology?

\subsection*{Acknowledgements}

We thank Tom Kennedy for useful conversations, and the referees for many useful suggestions regarding exposition and
inclusion of references.

\pdfbookmark[1]{References}{ref}
\LastPageEnding

\end{document}